\documentclass[11pt]{article}
\usepackage{}
\usepackage{amsmath}
\usepackage{amsfonts}
\usepackage{graphicx}
\usepackage{amssymb}

\allowdisplaybreaks

\newtheorem{condition**}{A*}
\newtheorem{condition***}{C*}
\newtheorem{condition*}{C}

\def\sX {{\mathcal X}}
\def\ro{{\rightarrow}}

\newtheorem{theorem}{Theorem}[section]

\newtheorem{lemma}[theorem]{Lemma}

\newtheorem{definition}[theorem]{Definition}

\newtheorem{remark}[theorem]{Remark}

\numberwithin{equation}{section}

\def\wt{\widetilde}
\def\qed{{\hfill $\Box$ \bigskip}}

\def\pf{\noindent{\bf Proof.} }

\begin{document}
\title{Reflected Backward Stochastic Differential Equation with Rank-based Data}
\author{Zhen-Qing Chen\thanks{Research supported in part by  Simons Foundation Grant 520542
and Victor Klee Faculty Fellowship at UW}
\quad \hbox{ and } \quad
  Xinwei Feng}

\date{July 11, 2020}

\maketitle
\renewcommand{\thefootnote}{\fnsymbol{footnote}}

\begin{abstract}
In this paper, we study reflected backward stochastic differential equation  (reflected BSDE  in abbreviation)
with rank-based data in a Markovian framework; that is, the solution
 to the reflected BSDE is
 above a  prescribed boundary process in a minimal fashion
and the generator and terminal value of the reflected BSDE depend on the solution
of another  stochastic differential equation  (SDE  in abbreviation)
 with rank-based drift and diffusion coefficients. We  derive regularity properties of the solution to such reflected BSDE, and show that the solution at the initial starting time $t$ and position $x$, which is a deterministic function, is the unique viscosity solution to some obstacle problem (or variational inequality) for the corresponding parabolic partial differential equation.
\end{abstract}



\section{Introduction}\label{S:1}

Non-linear backward stochastic differential equation was first investigated
by Pardoux and Peng \cite{PP1} in 1990. They obtained the existence and uniqueness of the solution under the  Lipschitz condition on the  generator function.  Reflected BSDE was studied in \cite{EKP},
where the solution is required to be above a certain prescribed process
and so an additional increasing process is needed in the equation to achieve this in a minimal fashion. Later, Cvitani\'{c} and Karatzas \cite{CK} studied  BSDE with two reflecting
barriers, that is, in addition to meet a target random data at the terminal time $T$, the solution
needs to remain between two processes.  For further work on BSDE, the reader is referred to El Karoui, Peng and Quenez \cite{EPQ}, Ma and Cvitani\'{c} \cite{MC}, Pardoux \cite{Pardoux}, Pardoux and Zhang \cite{PZ} and the references therein.

It has been widely recognized that BSDE provides a useful framework for formulating problems in many fields such as financial mathematics, stochastic optimal control and partial differential equation (PDE in abbreviation); see Chen and Epstein \cite{CE}, El Karoui, Peng and Quenez \cite{EPQ}, El Karoui and Quenez \cite{EQ}, Pardoux and Peng \cite{PP2}, Peng \cite{PENG} and the references therein. In \cite{CF}, we studied BSDE with the generator and terminal value depending on the solution of SDE with rank-based drift coefficients, and gave
  a probabilistic representation to the solution of the corresponding semi-linear backward parabolic PDE in cones with Neumann boundary condition.
 Application to  European option pricing problem with capital size based stock prices is also given. Paper \cite{CF} is motivated by the fact that the stock price of a company with large capital asset tends to move differently than that of a company with small capital asset.

 In rank-based SDEs,  drift and diffusion coefficients of each component are determined by its rank in  all components of the solution, which are piecewise constant. For SDEs with piecewise constant coefficients, Bass and Pardoux \cite{BP} proved the weak existence and uniqueness in law. Banner, Fernholz and Karatzas \cite{BFK} introduced an equity model that the growth rate and volatility of the stock depend on its rank in the market and studied the existence and uniqueness of the solution to the corresponding system of SDE. Shkolnikov \cite{Shkolnikov}  studied processes driven by independent identically distributed L\'{e}vy processes with rank-based coefficients and obtained the existence and uniqueness  of weak solution. For strong solution, Fernholz, Ichiba, Karatzas and Prokaj \cite{FIK} constructed a two dimensional process with rank-based coefficients. It is extended to finite or countably infinite system
in Ichiba, Karatzas and Shkolnikov \cite{IKS}, at least up until the first time that three particles collide simultaneously. The interested readers are referred to Ichiba and Karatzas \cite{IK}, Karatzas, Pal and Shkolnikov \cite{KPS}, Sarantsev \cite{Sarantsev1} and the references therein for more information about rank-based SDEs.

This paper is concerned with reflected BSDE with rank-based data, namely, the generator and terminal value of reflected BSDE depend on the solution from an SDE with rank-based drift and diffusion coefficients. First, we study the existence, uniqueness and regularity of solution of reflected BSDE with rank-based data. For this, we investigate regularity properties of the solution of rank-based SDE,  which is of interest in its own. Next, we show that the solution of the reflected BSDE provides a stochastic representation of the unique (under certain growth condition at infinity) viscosity solution to an obstacle problem for PDE in simplex with Neumann boundary conditions. Note that because of the dependence on the solution of rank-based SDE, the form of obstacle problem in this paper is defined on some simplex and is different from those in literature. Finally, a super-hedging problem for American option with capital size based stock prices is investigated which also serves as the motivation of the paper.

\medskip

Denote by $\Pi^n$ the open set $ \{ {x}\in\mathbb{R}^n: {x}_1 >  {x}_2>\cdots >  {x}_n\}$. Let
$$
  F_i:= \{ {x}\in \partial \Pi^n:  {x}_1 >  {x}_2
  > \cdots > {x}_i= {x}_{i+1}> \cdots >  {x}_n\}
 \quad \hbox{and} \quad \Gamma^n:=\Pi^n \cup (\cup_{i=1}^{n-1}F_i).
$$
Let $W$ be an $n$-dimensional Brownian motion,
and $\{X^{t, x}(s)=(X^{t, x}_1(s), \dots, X_n^{t, x}(s)) ; s\in [t, T]\}$  $\mathbb{R}^n$-valued
diffusion process
starting from $x$ at initial time $t$ with rank-based piecewise constant drifts $\delta_j$
and diffusion coefficients  $\sigma_j$ (see \S \ref{S:2.2} below for details).
The ranked process of $\{X^{t, x}(s); s\in [t, T]\}$ (from largest to the smallest) is
denoted by $\{\wt X^{t, x}(s); s\in [t, T]\}$.

Suppose $g:\Gamma^n\rightarrow\mathbb{R}$,
$F:[0,T]\times\Gamma^n\times\mathbb{R}\times\mathbb{R}^n\rightarrow\mathbb{R}$, and $h:[0,T]\times\Gamma^n\rightarrow\mathbb{R}$ are functions satisfying the following  conditions:

\begin{description}
\item(\textbf{H1}) $F(t, x, y, z)$ is jointly continuous on $[0,T]\times \Gamma ^n\times\mathbb{R}\times\mathbb{R}^n$,
and there exists a constant $c>0$ so that for any $t\in[0,T]$, $x\in\Gamma^n$, $y,y'\in\mathbb R$ and $z,z'\in\mathbb R^n$,
\begin{equation}\label{12}
\big|F(t,x,y,z)-F(t,x,y',z')\big|\leq c\big(|y-y'|+|z-z'|\big),
\end{equation}
and
\begin{equation}\label{13}
\big|F(t,x,0,0)\big|\leq c\big(1+|x|\big).
\end{equation}

\item(\textbf{H2})  There exists a constant $c$ so that
\begin{equation}\label{14}
\big|g(x)-g(x')\big|\leq c|x-x'|
\quad \hbox{for } x,x'\in\Gamma^n.
\end{equation}

\item(\textbf{H3}) $h:[0,T]\times\Gamma^n\rightarrow\mathbb{R}$ is jointly continuous and for some $c>0, p\in\mathbb{N}$ satisfies
\begin{equation}\label{21}h(t,x)\leq c\big(1+|x|^p\big),\ \ t\in[0,T],\ x\in\Gamma^n,
\end{equation}
and
\begin{equation}\label{22}h(T,x)\leq g(x),\ x\in\Gamma^n.
\end{equation}
\end{description}

For each $(t, x)\in [0, T]\times \Gamma^n$,
we consider the following reflected BSDE for $s\in [t, T]$:
\begin{equation} \label{e:1.1}
\left\{
\begin{aligned}
 Y^{t, {x}}(s) =& \,  \, g\big(\widetilde{X}^{t, {x}}(T)\big)+\int_s^T F\big(r, \widetilde X^{t, {x}}(r),Y^{t, {x}}(r), \bar Z^{t, {x}}(r)\big)dr +K^{t, {x}}(T)-K^{t, {x}}(s)   \\
&   \,     -\int_s^T Z^{t, {x}}(r)\cdot dW(r), \\
Y^{t, {x}}(s)   \geq & \, \, h\big(s, {X}^{t, {x}}(s)\big) \quad \hbox{for all } s\in [t, T],
\end{aligned}
\right.
\end{equation}
where
$\bar Z^{t, {x}}_j(r):=\sum_{i=1}^n Z^{t, {x}}_i(r)1_{\{ X^{t, {x}}_i(r)=\widetilde X^{t, {x}}_j(r)\}}$ and
$s\mapsto K^{t, x}(s)$ is a non-decreasing continuous process that increases only when
$Y^{t, x}(s)=h\big(s, {X}^{t, {x}}(s)\big)$. We emphasize that a solution to reflected BSDE \eqref{e:1.1} is a triple of progressively measurable processes
$\{(Y^{t, {x}}(s), Z^{t, {x}}(s), K^{t, {x}} (s)); s\in [t, T]\}$
taking values in $\mathbb{R} \times \mathbb{R}^n \times \mathbb{R}^+$ so that equation
\eqref{e:1.1} holds.
For BSDE \eqref{e:1.1} related to option pricing, $Z^{t, x}_i(r)$ corresponds to number of shares invested in the $i$th stock at time $r$, and so $\bar Z^{t, x}_j(r)$ is the number of shares invested in the $j$th ranked stock at time $r$.

Under conditions ({\bf H1})-({\bf H3)} and that the sequence $\{\sigma_1^2, \dots, \sigma_n^2\}$ is concave (see Definition \ref{D:2.1}(i) for its definition),
 for every $(t, x)$, we obtain
in Theorem \ref{th6} that
the above reflected BSDE has a unique strong solution $(Y^{t, {x}}, Z^{t, {x}}, K^{t, {x}})$.
 By Kolmogorov's 0-1 law, $u(t, x):= Y^{t, x}(t)$ is a deterministic number and hence a function of $(t, x)$. The following existence and uniqueness of viscosity solution for \eqref{e:1.7}
are the other two main results of this paper.

\begin{theorem}\label{th8}
Suppose that {\bf(H1)}-{\bf (H3)} hold and that the sequence $\{\sigma_1^2, \dots, \sigma_n^2\}$ is concave. The function $u(t,x)= Y^{t, x}(t)$ is a viscosity solution of the following obstacle problem:
\begin{equation}\label{e:1.7}
\left\{
\begin{aligned}
&\bigg(u(t, {x})-h(t, {x})\bigg)\wedge\bigg( -\frac{\partial u}{\partial t}(t, {x})-\mathcal{L} u(t, {x})-F\big(t, {x},u(t, {x}),\sigma(\nabla u)(t, {x})\big)\bigg)=0 \\
& \hskip 3.5 truein  \hbox{for }  t\in[0,T) \hbox{ and }  {x}\in\Pi^n , \\
&u(T, {x})=g( {x}) \qquad \hbox{for }  {x}\in\Gamma^n,\\
&\frac{\partial u}{\partial  {x}_{i+1}}(t, {x})=\frac{\partial u}{\partial  {x}_i}(t, {x}) \qquad \hbox{for } t\in[0,T) \hbox{ and }
 {x}\in F_i,\ i=1,\ldots,n-1.
\end{aligned}
\right.
\end{equation}
where
\begin{equation}\nonumber
\mathcal{L}=\frac{1}{2}\sum_{i=1}^n\sigma_i^2\frac{\partial^2}{\partial x_{i}^2}+\sum_{i=1}^n\delta_i\frac{\partial}{\partial x_{i}} .
\end{equation}
\end{theorem}

\begin{theorem}\label{th9}
Suppose {\bf(H1)}-{\bf(H3)} hold. If for each $R>0$, there exists a positive function $\eta_R$ on $[0, \infty)$ with $\lim_{r\ro 0^+} \eta_R (r) =0$ such that
\begin{equation}\label{27}
\big|F(t,x,y,z)-F(t,x',y,z)\big|\leq\eta_R\big(|x-x'|(1+|z|)\big),
\end{equation}
for all $t\in[0,T]$, $z\in \mathbb{R}^n$, $x, x'\in \Gamma^n$ and $y\in \mathbb{R}$ having
$\max\{|x|,|x'|,|y|\}\leq R$,
 then there exists at most one viscosity solution $u$ of \eqref{e:1.7} such that
\begin{equation}\label{28}
\lim_{|x|\rightarrow \infty}\big|u(t,x)\big|e^{-A\log^2|x|}=0,
\end{equation}
uniformly in $t\in[0,T]$, for some $A>0$.
\end{theorem}

 In the process, we establish regularity results for the solution of rank-based SDEs. Our approach (see Theorem \ref{th1} and Lemma \ref{lemma1}) may be useful to study other SDEs with non-Lipschitz continuous coefficients. A key step in proving the uniqueness of viscosity solution to the obstacle problem \eqref{e:1.7} is Lemma \ref{lemma2}, which gives a characterization of the difference between a viscosity subsolution and a viscosity supersolution of \eqref{e:1.7}.

\medskip

The rest of this paper is organized as follows. In Section \ref{S:2}, we introduce the notations and preliminaries on rank-based SDE and BSDE. In Section
\ref{S:3}, we study reflected BSDE with rank-based data. Existence, uniqueness and   probabilistic representation of the solution of the obstacle problem \eqref{e:1.7}
are established in Section \ref{S:4}.
In Section \ref{S:5}, we study American option pricing where the stock prices  have rank-based drift and volatility coefficients.

\section{Preliminaries}\label{S:2}

\subsection{Notations}\label{S:2.1}
Let $n\geq1$ be the dimension.
 Denote by $\partial(\Pi^n)$ the boundary of $\Pi^n$, and $S(n)$ the set of $n\times n$ symmetric matrices. For a vector or matrix $a$, $a^{tr}$ stands for its transpose. The Euclidean inner product between two vectors $x, y\in \mathbb{R}^n$ will be denoted by $x\cdot y$. Let $(\Omega,\mathcal{F},\mathbb{P})$ be a   complete probability space  on which defines
a standard $n$-dimensional Brownian motion
 $W=(W_1,W_2,\cdots,W_n)$. Denote by $\{\mathcal{F}^W_t\}$ the minimal augmented filtration generated by $\{W(t)\}$
 such that $\mathcal{F}^W_0$ contains all the $\mathbb{P}$-null subsets of $\mathcal{F}^W_\infty$. Fix $T>0$. For $n\in\mathbb{N}$ and $p\geq1$, consider the following spaces of random variables or stochastic processes:

\begin{itemize}
  \item $L^p(\mathcal{F}^W_T;\mathbb{R}^n)$:
the set of $\mathbb{R}^n$-valued $\mathcal{F}^W_T$-measurable random variables $\xi$ with
$\mathbb{E}\big[|\xi|^p\big]< \infty$.
  \item $M^p\big([0,T];\mathbb{R}^n\big)$: the set of $n$-dimensional progressively measurable processes $\{\xi_t,0\leq t\leq T\}$ with $\mathbb{E}\big[\int_0^T|\xi_t|^pdt\big]<\infty.$
  \item $S^p\big([0,T];\mathbb{R}^n\big)$: the set of $n$-dimensional progressively measurable processes $\{\xi_t,0\leq t\leq T\}$ with $\mathbb{E}\big[\sup_{0\leq t\leq T}|\xi_t|^p\big]<\infty.$
\end{itemize}

Denote by $C\big([0,T]\times\mathbb{R}^n\big)$ the space of continuous functions on $[0, T] \times\mathbb{R}^n$, and $C^{1,2}\big([0,T)\times\Gamma^n\big)$ the space of functions $u(t,x):[0,T)\times\Gamma^n\rightarrow\mathbb{R}$ that are  continuously differentiable in $t$ and twicely continuously differentiable in $x$. For $u\in C\big([0,T]\times\mathbb{R}^n\big)$,  $\bar{D}^{2,+}_u(t,x)$ is the parabolic superjet of $u$ at $(t,x)$; that is, $\bar{D}^{2,+}_u(t,x)$ is the set of triples $(p,q,X)\in\mathbb{R}\times\mathbb{R}^n\times S(n)$ so that
$$
u(s,y)\leq u(t,x)+p(s-t)+ q \cdot(y-x)+\frac{1}{2}(y-x)^{tr}\cdot X(y-x)+o\big(|s-t|+|y-x|^2\big)
$$
for $(s,y)\in[0,T]\times\mathbb{R}^n$.
Here the notation $o (\delta)$ means a quantity $f(\delta)$ such that $\lim_{\delta \to 0} f(\delta)/\delta =0$.
Similarly,  $\bar{D}^{2,-}_u(t,x)$ denotes the parabolic subjet of $u$ at $(t,x)$; that is, $\bar{D}^{2,-}_u(t,x)$ is the set of triples $(p,q,X)\in\mathbb{R}\times\mathbb{R}^n\times S(n)$ such that
$$
u(s,y)\geq u(t,x)+p(s-t)+ q \cdot (y-x)+\frac{1}{2}(y-x)^{tr}\cdot X(y-x)+o\big(|s-t|+|y-x|^2\big)
$$
for $(s,y)\in[0,T]\times\mathbb{R}^n$.
For $0\leq t<T$, denote by $C\big([t,T];\mathbb{R}^n\big)$ the family of $\mathbb{R}^n$-valued continuous functions  defined on $[t,T]$.

We will use $C$ to denote a constant whose value may change from line to line. For $a, b\in \mathbb{R}$, $a\vee b:=\max \{a, b\}$
and $a\wedge b  := \min\{a, b\}$. We use $\langle X, Y\rangle$ to denote the
quadratic covariation process of two continuous semimartingales $X$ and $Y$.

\subsection{Rank-based SDE}\label{S:2.2}

In this subsection, we present some properties of SDE with rank-based drift and diffusion coefficients. For any initial value $(t,x)\in[0,T]\times\Gamma^n$ and $i=1,\cdots,n$, consider the following SDE
\begin{equation}\label{1}
X_i^{t,x}(s)=x_i+\int_t^{s\vee t}\sum_{j=1}^n\textbf{1}_{\{X_i^{t,x}(r)=X_{(j)}^{t,x}(r)\}}\delta_jdr+\int_t^{s\vee t}\sum_{j=1}^n\textbf{1}_{\{X_i^{t,x}(r)=X_{(j)}^{t,x}(r)\}}\sigma_jdW_i(r),
\end{equation}
 where  $\delta_j,\ j=1,\cdots,n,$ are real constants; $\sigma_j,\ j=1,\cdots,n,$ are strictly positive real constants and
\begin{equation}\label{4}
X_{(1)}^{t,x}(r)\geq X_{(2)}^{t,x}(r)\geq\cdots\geq X_{(n)}^{t,x}(r)
\end{equation}
are the ordered configuration of $\{ X_1^{t,x}(r), X_2^{t,x}(r),
\dots ,  X_n^{t,x}(r)\}$
at time $r$, with ties resolved by resorting to the lowest index. For instance, we set $$X_{(i)}^{t,x}(r)=X_i^{t,x}(r),\ i=1,\cdots,n, \quad  \text{whenever}\ X_1^{t,x}(r)=\cdots=X_n^{t,x}(r).$$
The processes $X_i^{t,x}$, $i=1,\cdots,n,$ are called named particles and $X_{(j)}^{t,x}$, $j=1,\cdots,n,$ are called ranked particles.

\begin{definition}\label{D:2.1}
 \begin{description}
\item{\rm (i)} A finite sequence $\{a_1,\ldots,a_n\}$ of real numbers is called concave, if for every three consecutive elements $a_i,a_{i+1},a_{i+2}$, we have $a_{i+1}\geq\frac{1}{2}(a_i+a_{i+2})$, $i=1,\ldots,n-2$.

\item{\rm (ii)} A triple collision at time $t$ occurs if there exists some rank $j=2,\ldots,n-1$, so that $X_{(j-1)}(r)=X_{(j)}(r)=X_{(j+1)}(r).$
 \end{description}
\end{definition}

\begin{theorem}\emph{(\cite[Theorem 2]{IKS} and \cite[Theorem 1.4]{Sarantsev2})} \label{T:2.2}
Suppose the sequence $\{\sigma_1^2,\ldots,\sigma_n^2\}$ is concave.
Then with probability one, there are no triple collisions at any time $t > 0$ and there exists a unique strong solution of the system \eqref{1} defined for all $t\geq0$.
\end{theorem}
It follows from \cite{BG} (see also \cite[Section 2]{IKS}) that the ranked particles have the following representation: for $j=1,\cdots,n,$
$$dX_{(j)}^{t,x}(s)=\delta_jds+ \sigma_jd\beta_j(s)+\frac{1}{2}\big(d\Lambda^{j,j+1}(s)-d\Lambda^{j-1,j}(s)\big),\ \ s\geq t,$$
where $\Lambda^{j,j+1}(s)$, $j=1,\cdots,n-1,$ denotes the local times accumulated at the origin by the nonnegative semimartingales (called gap processes)
$$
G_j^{t,x}(\cdot)=X_{(j)}^{t,x}(\cdot)-X_{(j+1)}^{t,x}(\cdot),\quad \ j=1,\cdots,n-1,
$$
over the time interval $[0,s]$ and
\begin{equation}\label{e:2.3}
\beta_j(\cdot):=\sum_{i=1}^n\int_0^{\cdot}\textbf{1}_{\{X_i^{t,x}(r)=X_{(j)}^{t,x}(r)\}}dW_i(r),\quad j=1,\cdots,n.
\end{equation}
are independent Brownian motions.

The process of ranked particles is a semimartingale (normally) reflected Brownian motion in the open set $\Gamma^n$. That is, the process $\big(X_{(1)}^{t,x}(\cdot),\cdots,X_{(n)}^{t,x}(\cdot)\big)$ behaves like an $n$-dimensional Brownian motion with constant drift and covariance matrix in the open set $\Pi^n$ and is normally reflected on the faces $F_i$, $i=1,\cdots,n-1$. The gap process $\big(G_1^{t,x}(\cdot),\cdots,G_{n-1}^{t,x}(\cdot)\big)$ is an obliquely reflected Brownian motion in the $(n-1)$-dimensional non-negative orthant $(\mathbb R^+)^{n-1}$.
We refer the interested readers to Dai and Williams \cite{DW} for the properties of semimartingale reflected Brownian motion in the convex  polyhedrons and Harrison and Reiman \cite{HR}, Williams \cite{Williams} for properties of reflected Brownian motions in orthants.

Next we study regularity properties of the solution of SDE \eqref{1}, which are needed in the study of that of reflected BSDE with rank-based data. To the best of our knowledge, this is the first result on regularity of solutions of \eqref{1} and the method may be applicable to other SDE with non-Lipschitz continuous drift and diffusion coefficients.

\begin{theorem}\label{th1}
Under the conditions of Theorem \ref{T:2.2}, for all $p\geq1$, there exists a constant $C>0$ depending on $\big(p,T,\{\delta_j\},\{\sigma_j\}\big)$ such that for any $x,\{{x^{(k)}}\}_{k=1}^\infty\in\Gamma^n$ and any $t,\{{t^{(k)}}\}_{k=1}^\infty\in[0,T]$ with $\lim_{k\rightarrow\infty}(|{t^{(k)}}-t|+|{x^{(k)}}-x|)=0$, we have for $i=1,\cdots,n,$
\begin{equation}\label{2}
\mathbb{E}\bigg[\sup_{0\leq s\leq T}\big|X_i^{t,x}(s)\big|^p\bigg]\leq C\big(1+|x|^p\big),
\end{equation}
and
\begin{equation}\label{3}
\lim_{k\rightarrow\infty}\mathbb{E}\bigg[\sup_{0\leq s\leq T}\big|X_i^{t,x}(s)-X_i^{{t^{(k)}},{x^{(k)}}}(s)\big|^p\bigg]=0.
\end{equation}
\end{theorem}

To prove Theorem \ref{th1}, we  need the following lemma.

\begin{lemma}\label{lemma1}
Let $b(x)$ be a bounded $\mathbb{R}^n$-valued Borel measurable function defined on $\mathbb{R}^n$ and $\sigma(x)$ a bounded $\mathbb{R}^{n\times n}$-valued Borel measurable function defined on $\mathbb{R}^n$ that is uniformly elliptic, i.e., there exists a constant $\alpha>0$ so that
$$
y \cdot (\sigma(x)y) \geq\alpha|y|^2 \quad \hbox{for every } x, y \in \mathbb{R}^n.
$$
Fix $T>0$ and $(t,x)\in[0,T]\times\mathbb{R}^n$.  Suppose  SDE
\begin{equation}\label{5}
X^{t,x}(s)=x+\int_t^{s\vee t}b\big(X^{t,x}(u)\big)du+\int_t^{s\vee t}\sigma\big(X^{t,x}(u)\big)dW(u)
\end{equation}
has a strong solution and pathwise uniqueness holds.
Let $\{ b_m (x); m\geq 1\}$ be a sequence of $\mathbb{R}^n$-valued uniformly bounded Lipschitz continuous functions that converges to $b(x)$ almost everywhere on $\mathbb{R}^n$ with respect to the Lebesgue measure,  and $\{\sigma_m(x); m\geq 1\}$ a sequence of $\mathbb{R}^{n\times n}$-valued uniformly bounded Lipschitz continuous functions so that for some $\alpha >0$,
$$
y\cdot (\sigma_m(x)y)
\geq\alpha|y|^2 \quad \hbox{for all } x, y\in \mathbb{R}^n \hbox{ and }
m\geq 1,
$$
and $\lim_{m\rightarrow\infty}\sigma_m(x)=\sigma(x)$ almost everywhere on $\mathbb{R}^n$ with respect to the Lebesgue measure.
Denote by $X^{t,x,m}$ the unique strong solution of \eqref{5} but with $b_m$ in place of $b$ and $\sigma_m$ in place of $\sigma$.
Then for every $q\geq 1$ and any compact subset $\mathcal{K}$ of $\mathbb{R}^n$,
\begin{equation}\label{6}
\lim_{m\rightarrow \infty}\sup_{0\leq t\leq T}\sup_{x\in\mathcal{K}}\mathbb{E}\bigg[\sup_{0\leq s\leq T}\big|X^{t,x,m}(s)-X^{t,x}(s)\big|^q\bigg]=0 .
\end{equation}
\end{lemma}
\pf
 The proof adopts some idea from Kaneko and Nakao \cite{KN}, with necessary modifications. Suppose the conclusion of the lemma is not true. Then there exists a positive constant $\varepsilon$ and a subsequence of $\{m\}$ (still denoted by $\{m\}$), a sequence $\{t^{(m)}\}$ contained in $[0,T]$ and a sequence $\{x^{(m)}\}$ contained in some compact subset $\mathcal K$ of $\mathbb{R}^n$ such that
$$\inf_m\mathbb{E}\bigg[\sup_{0\leq s\leq T}\big|X^{t^{(m)},x^{(m)},m}(s)-X^{t^{(m)},x^{(m)}}(s)\big|^q\bigg]\geq\varepsilon.$$
Without loss of generality, we may assume that $\{t^{(m)}\}$ converges to $t$ in $[0,T]$ and $\{x^{(m)}\}$ converges to $x$ in $\mathbb{R}^n$.

First, since $b$, $b_m$, $\sigma$ and $\sigma_m$ are uniformly bounded,  there exists a constant $C$ depending on $(T,\mathcal{K})$ such that for $0\leq r_1<r_2\leq T$,
$$
\sup_{0\leq t\leq T}\sup_{x\in\mathcal{K}}\mathbb{E}\bigg[\sup_{r_1\leq u_1,u_2\leq r_2}\big|X^{t,x}(u_2)-X^{t,x}(u_1)\big|^4\bigg]\leq C|r_2-r_1|^2,
$$
and
$$
\sup_m\sup_{0\leq t\leq T}\sup_{x\in\mathcal{K}}\mathbb{E}\bigg[\sup_{r_1\leq u_1,u_2\leq r_2}\big|X^{t,x,m}(u_2)-X^{t,x,m}(u_1)\big|^4\bigg]\leq C|r_2-r_1|^2.
$$
Thus the family of the processes $\big\{X^{t^{(m)},x^{(m)}}(s),X^{t^{(m)},x^{(m)},m}(s),W(s)\big\}_{m=1}^{ \infty},s\in[0,T]$ is tight (see, e.g., \cite[Theorems 1.4.2 and 1.4.3]{IW}).
By Skorohod embedding theorem, there exist some probability space $(\hat{\Omega},\hat{\mathcal{F}},\hat{\mathbb{P}})$ and a sequence of continuous stochastic processes $\big\{\hat{X}^{(m)}(s),\hat{Y}^{(m)}(s),\hat{W}^{(m)}(s)\big\}_{m=0}^{ \infty}$ on it which enjoy the following properties:
\begin{description}
  \item[(i)] The probability law of $\{\hat{X}^{(m)},\hat{Y}^{(m)},\hat{W}^{(m)}\big\}$ coincides with the law of $\{X^{t^{(m)},x^{(m)}},X^{t^{(m)},x^{(m)},m},W\}$ for each $m=1,2,\ldots$
  \item[(ii)] There exists a subsequence $\{m_j\}$ such that $\{\hat{X}^{(m_j)},\hat{Y}^{(m_j)},\hat{W}^{(m_j)}\big\}$ converges to $\{\hat{X}_0,\hat{Y}_0,\hat{W}_0\big\}$ uniformly on every finite time interval a.s.
\end{description}
For simplicity, we write $t^{(m_j)}=t^{(m)}$, $x^{(m_j)}=x^{(m)}$, $\hat{X}^{(m_j)}=\hat{X}^{(m)}$, $\hat{Y}^{(m_j)}=\hat{Y}^{(m)}$ and $\hat{W}^{(m_j)}=\hat{W}^{(m)}$. By virtue of uniformly integrability, we obtain
\begin{equation}\label{7}
\begin{aligned}
\varepsilon&\leq\liminf_{m\rightarrow \infty}\mathbb{E}\bigg[\sup_{0\leq s\leq T}\big|X^{t^{(m)},x^{(m)},m}(s)-X^{t^{(m)},x^{(m)}}(s)\big|^q\bigg]\\
&=\liminf_{m\rightarrow \infty}\hat{\mathbb{E}}\bigg[\sup_{0\leq s\leq T}\big|\hat{X}^{(m)}(s)-\hat{Y}^{(m)}(s)\big|^q\bigg]\\
&=\hat{\mathbb{E}}\bigg[\sup_{0\leq s\leq T}\big|\hat{X}_0(s)-\hat{Y}_0(s)\big|^q\bigg].\\
\end{aligned}
\end{equation}
On the other hand, because of the coincidence (i) of probability law, we have for $m=1,2,\ldots$,
$$\hat{X}^{(m)}(s)=x^{(m)}+\int_{t^{(m)}}^{s\vee t^{(m)}}b\big(\hat{X}^{(m)}(u)\big)du+\int_{t^{(m)}}^{s\vee t^{(m)}}\sigma\big(\hat{X}^{(m)}(u)\big)d\hat{W}^{(m)}(u),$$
and $$\hat{Y}^{(m)}(s)=x^{(m)}+\int_{t^{(m)}}^{s\vee t^{(m)}}b_m\big(\hat{Y}^{(m)}(u)\big)du+\int_{t^{(m)}}^{s\vee t^{(m)}}\sigma_m\big(\hat{Y}^{(m)}(u)\big)d\hat{W}^{(m)}(u).$$
Now we are going to take the limit $m\rightarrow \infty$ in the above two equations. First, we will show how to take the limit in the drift term.

 If $s=t$, then $\int_t^{s\vee t}b\big(\hat{X}_0(u)\big)du=0$. Since $$\lim_{m\rightarrow \infty}\big|\int_{t^{(m)}}^{s\vee t^{(m)}}b\big(\hat{X}^{(m)}(u)\big)du\big|
 \leq C\lim_{m\rightarrow \infty}|t\vee t^{(m)}-t^{(m)}|=0,$$ we have,
 \begin{equation}\label{lemma2.4-1}
\lim_{m\rightarrow \infty}\hat{\mathbb{E}}\left[\bigg|\int_{t^{(m)}}^{s\vee t^{(m)}}b\big(\hat{X}^{(m)}(u)\big)du-\int_t^{s\vee t}b\big(\hat{X}_0(u)\big)du\bigg|\right]=0.
\end{equation}

 If $s< t$, then $\int_t^{s\vee t}b\big(\hat{X}_0(u)\big)du=0$. Since $\lim_{m\rightarrow \infty}t^{(m)}=t$, there exists a sufficiently large $M\in\mathbb{N}$ such that for all $ m\geq M$, $s\leq t^{(m)}.$ Thus, $\int_{t^{(m)}}^{s\vee t^{(m)}}b\big(\hat{X}^{(m)}(u)\big)du=0.$
Therefore, for $s<t$, we have \eqref{lemma2.4-1}.

 If $s> t$, from $\lim_{m\rightarrow \infty}t^{(m)}=t$, we know that there exists a sufficiently large $M\in\mathbb{N}$ such that for all $ m\geq M$, $s\geq t^{(m)}.$ Fix $m_0$, we have

\begin{eqnarray*}
& &\hat{\mathbb{E}}\left[\bigg|\int_{t^{(m)}}^sb\big(\hat{X}^{(m)}(u)\big)du-\int_t^sb\big(\hat{X}_0(u)\big)du\bigg|\right]\\
&\leq&\hat{\mathbb{E}}\left[\bigg|\int_0^s\big[\textbf{1}_{[t^{(m)},s]}(u)b\big(\hat{X}^{(m)}(u)\big)-\textbf{1}_{[t^{(m)},s]}(u)
b_{m_0}\big(\hat{X}^{(m)}(u)\big)\big]du\bigg|\right]\\
& &+\hat{\mathbb{E}}\left[\bigg|\int_0^s\big[\textbf{1}_{[t^{(m)},s]}(u)b_{m_0}\big(\hat{X}^{(m)}(u)\big)-\textbf{1}_{[t^{(m)},s]}(u)b_{m_0}\big(\hat{X}_0(u)\big)\big]du\bigg|\right]\\
& &+\hat{\mathbb{E}}\left[\bigg|\int_0^s\big[\textbf{1}_{[t^{(m)},s]}(u)b_{m_0}\big(\hat{X}_0(u)\big)-\textbf{1}_{[t,s]}(u)b_{m_0}\big(\hat{X}_0(u)\big)\big]du\bigg|\right]\\
& &+\hat{\mathbb{E}}\left[\bigg|\int_0^s\big[\textbf{1}_{[t,s]}(u)b_{m_0}\big(\hat{X}_0(u)\big)-\textbf{1}_{[t,s]}(u)b\big(\hat{X}_0(u)\big)\big]du\bigg|\right]\\
&\leq&\hat{\mathbb{E}}\left[\int_{t^{(m)}}^s\big|b\big(\hat{X}^{(m)}(u)\big)-b_{m_0}\big(\hat{X}^{(m)}(u)\big)\big|du\right]
+\hat{\mathbb{E}}\left[\int_0^s\big|b_{m_0}\big(\hat{X}^{(m)}(u)\big)-b_{m_0}\big(\hat{X}_0(u)\big)\big|du\right]\\
& &+\hat{\mathbb{E}}\left[\int_0^s\big|b_{m_0}\big(\hat{X}_0(u)\big)\big|\big|\textbf{1}_{[t^{(m)},s]}(u)-\textbf{1}_{[t,s]}(u)\big|du\right]
+\hat{\mathbb{E}}\left[\int_t^s\big|b_{m_0}\big(\hat{X}_0(u)\big)-b\big(\hat{X}_0(u)\big)\big|du\right]\\
&=:&I_1+I_2+I_3+I_4.
\end{eqnarray*}
Let $w(x)$ be a $[0,1]$-valued continuous function defined on $\mathbb{R}^n$ such that $w(0)=1$ and $w(x)=0$ for $|x|^2\geq1$.
Then for $R>0$, it follows from \cite[Theorem 2.2.4]{Krylov} that
\begin{eqnarray*}
I_1
&\leq&C\hat{\mathbb{E}}\left[\int_{t^{(m)}}^s\big[1-w\big(\frac{\hat{X}^{(m)}(u)}{R}\big)\big]du\right]+\hat{\mathbb{E}}
\left[\int_{t^{(m)}}^sw\big(\frac{\hat{X}^{(m)}(u)}{R}\big)
\big|b\big(\hat{X}^{(m)}(u)\big)-b_{m_0}\big(\hat{X}^{(m)}(u)\big)\big|du\right]\\
&\leq&C\hat{\mathbb{E}}\left[\int_{t^{(m)}}^s\big[1-w\big(\frac{\hat{X}^{(m)}(u)}{R}\big)\big]du\right]+C\bigg(\int_{B(0,R)}\big|b(y)-b_{m_0}(y)\big|^{n+1}dy\bigg)^{\frac{1}{n+1}},\\
\end{eqnarray*}
where $B(0,R)$ is the ball with center 0 and radius $R$ in $\mathbb{R}^n$. Therefore,
$$
\lim_{m\rightarrow \infty}I_1
\leq C\hat{\mathbb{E}}\left[\int_0^s\big[1-w\big(\frac{\hat{X}_0(u)}{R}\big)\big]du\right]+C\bigg(\int_{B(0,R)}\big|b(y)-b_{m_0}(y)\big|^{n+1}dy\bigg)^{\frac{1}{n+1}}.
$$
First let $m_0$ tends to $ \infty$ and then let $R\rightarrow\infty$, we have
$$\lim_{m\rightarrow \infty}I_1=0.$$
Similarly, we have
$$\lim_{m\rightarrow \infty}I_4=0.$$
From the convergence of $\hat{X}^{(m)}$ to $\hat{X}_0$, the continuity of $b_{m_0}$ and dominated convergence theorem, we have
$$\lim_{m\rightarrow \infty}I_2=0.$$
Finally,
$$\lim_{m\rightarrow \infty}I_3\leq\lim_{m\rightarrow \infty}C|t^{(m)}-t|=0.$$
Thus, if $s>t$, we have \eqref{lemma2.4-1}.
Therefore, for $0\leq s\leq T$,
$$\lim_{m\rightarrow \infty}\hat{\mathbb{E}}\left[\bigg|\int_{t^{(m)}}^{s\vee t{(m)}}b\big(\hat{X}^{(m)}(u)\big)du-\int_t^{s\vee t}b\big(\hat{X}_0(u)\big)du\bigg|\right]=0.$$
Similarly, we have
$$\lim_{m\rightarrow \infty}\hat{\mathbb{E}}\left[\bigg|\int_{t^{(m)}}^{s\vee t{(m)}}b_m\big(\hat{Y}^{(m)}(u)\big)du-\int_t^{s\vee t}b\big(\hat{Y}_0(u)\big)du\bigg|\right]=0.$$
Next, we consider the limit of diffusion term. The case $s\leq t$ is similar to that of drift term so is omitted here.
If $s> t$, from $\lim_{m\rightarrow \infty}t^{(m)}=t$, we know that there exists a sufficiently large $M\in\mathbb{N}$ such that for all $ m\geq M$, $s\geq t^{(m)}.$ Fix $m_0$, we have

\begin{eqnarray*}
& &\hat{\mathbb{E}}\left[\bigg|\int_{t^{(m)}}^s\sigma\big(\hat{X}^{(m)}(u)\big)d\hat{W}^{(m)}(u)
-\int_t^s\sigma\big(\hat{X}_0(u)\big)d\hat{W}_0(u)\bigg|^2\right]\\
&\leq&C\hat{\mathbb{E}}\left[\bigg|\int_0^s\big[\textbf{1}_{[t^{(m)},s]}(u)\sigma\big(\hat{X}^{(m)}(u)\big)-\textbf{1}_{[t^{(m)},s]}(u)
\sigma_{m_0}\big(\hat{X}^{(m)}(u)\big)\big]d\hat{W}^{(m)}(u)\bigg|^2\right]\\
& &+C\hat{\mathbb{E}}\left[\bigg|\int_0^s\textbf{1}_{[t^{(m)},s]}(u)\sigma_{m_0}\big(\hat{X}^{(m)}(u)\big)d\hat{W}^{(m)}(u)
-\int_0^s\textbf{1}_{[t,s]}(u)\sigma_{m_0}\big(\hat{X}^{(m)}(u)\big)d\hat{W}^{(m)}(u)\bigg|^2\right]\\
& &+C\hat{\mathbb{E}}\left[\bigg|\int_0^s\textbf{1}_{[t,s]}(u)\sigma_{m_0}\big(\hat{X}^{(m)}(u)\big)d\hat{W}^{(m)}(u)
-\int_0^s\textbf{1}_{[t,s]}(u)\sigma_{m_0}\big(\hat{X}_0(u)\big)d\hat{W}_0(u)\bigg|^2\right]\\
 & &+C\hat{\mathbb{E}}\left[\bigg|\int_0^s\big[\textbf{1}_{[t,s]}(u)\sigma_{m_0}\big(\hat{X}_0(u)\big)
-\textbf{1}_{[t,s]}(u)\sigma\big(\hat{X}_0(u)\big)\big]d\hat{W}_0(u)\bigg|^2\right]\\
&\leq&C\hat{\mathbb{E}}\left[\int_{t^{(m)}}^s\big|\sigma\big(\hat{X}^{(m)}(u)\big)
-\sigma_{m_0}\big(\hat{X}^{(m)}(u)\big)\big|^2du\right]\\
& &+C\hat{\mathbb{E}}\left[\int_0^s\big|\textbf{1}_{[t^{(m)},s]}(u)\sigma_{m_0}\big(\hat{X}^{(m)}(u)\big)
-\textbf{1}_{[t,s]}(u)\sigma_{m_0}\big(\hat{X}^{(m)}(u)\big)\big|^2du\right]\\
& &+C\hat{\mathbb{E}}\left[\bigg|\int_0^s\textbf{1}_{[t,s]}(u)\sigma_{m_0}\big(\hat{X}_m(u)\big)d\hat{W}^{(m)}(u)
-\int_0^s\textbf{1}_{[t,s]}(u)\sigma_{m_0}\big(\hat{X}_0(u)\big)\big]d\hat{W}_0(u)\bigg|^2\right]\\
& &+C\hat{\mathbb{E}}\left[\int_t^s\big|\sigma_{m_0}\big(\hat{X}_0(u)\big)-\sigma\big(\hat{X}_0(u)\big)\big|^2du\right]\\
&=:&J_1+J_2+J_3+J_4.
\end{eqnarray*}
Similar as $I_1$ and $I_4$, we have $\lim_{m\rightarrow \infty}J_1=\lim_{m\rightarrow \infty}J_4=0.$ It is easy to see that $\lim_{m\rightarrow \infty}J_2=0.$ Due to \cite[Lemma 5.2]{GM}, we have that $$\int_0^s\textbf{1}_{[t,s]}(u)\sigma_{m_0}\big(\hat{X}^{(m)}(u)\big)d\hat{W}^{(m)}(u)\rightarrow\int_0^s\textbf{1}_{[t,s]}(u)\sigma_{m_0}\big(\hat{X}_0(u)\big)d\hat{W}_0(u)$$in probability, then by uniformly integrable, we have $ \lim_{m\rightarrow \infty}J_3=0$.
Therefore, we have
$$\lim_{m\rightarrow \infty}\hat{\mathbb{E}}\left[\bigg|\int_{t^{(m)}}^{s\vee t^{(m)}}\sigma\big(\hat{X}^{(m)}(u)\big)d\hat{W}^{(m)}(u)-\int_t^{s\vee t}\sigma\big(\hat{X}_0(u)\big)d\hat{W}_0(u)\bigg|^2\right]=0.$$
Similarly,
$$\lim_{m\rightarrow \infty}\hat{\mathbb{E}}\left[\bigg|\int_{t^{(m)}}^{s\vee t^{(m)}}\sigma_m\big(\hat{Y}^{(m)}(u)\big)d\hat{W}^{(m)}(u)-\int_t^{s\vee t}\sigma\big(\hat{Y}_0(u)\big)d\hat{W}_0(u)\bigg|^2\right]=0.$$
Therefore, we have that both $\hat{X}_0$ and $\hat{Y}_0$ are solutions of
$$X(s)=x+\int_t^{s\vee t}b\big(X(u)\big)du+\int_t^{s\vee t}\sigma\big(X(u)\big)d\hat{W}_0(u).$$
From the pathwise uniqueness of solutions for the above SDE, we have $\hat{X}_0(s)=\hat{Y}_0(s)$ almost surely for $0\leq s\leq T.$ This contradicts \eqref{7} and the proof is complete. \qed
\medskip

\noindent{\bf \textbf{Proof of Theorem \ref{th1}:}}
For any $i=1,\cdots,n$,
\begin{eqnarray*}
\big|X_i^{t,x}(s)\big|^p&=&\bigg|x_i+\int_t^{s\vee t}\sum_{j=1}^n\textbf{1}_{\{X_i^{t,x}(r)=X_{(j)}^{t,x}(r)\}}\delta_j dr+ \int_t^{s\vee t}\sum_{j=1}^n\textbf{1}_{\{X_i^{t,x}(r)=X_{(j)}^{t,x}(r)\}}\sigma_jdW_i(r)\bigg|^p\\
&\leq&C\bigg(|x|^p+|s\vee t|^p+\big|\int_t^{s\vee t}\sum_{j=1}^n\textbf{1}_{\{X_i^{t,x}(r)=X_{(j)}^{t,x}(r)\}}\sigma_j dW_i(r)\big|^p\bigg).\\
\end{eqnarray*}
By Burkholder-Davis-Gundy inequality, we have
\begin{eqnarray*}
\mathbb{E}\bigg[\sup_{0\leq s\leq T}\big|X_i^{t,x}(s)\big|^p\bigg]&\leq&C\bigg(|x|^p+|T|^p+\mathbb{E}\big[\sup_{0\leq s\leq T}\big|\int_t^{s\vee t}\sum_{j=1}^n\textbf{1}_{\{X_i^{t,x}(r)=X_{(j)}^{t,x}(r)\}}\sigma_jdW_i(r)\big|^p\big]\bigg).\\
&\leq& C\big(1+|x|^p\big).\
\end{eqnarray*}
For $i=1,\ldots,n$ and $x\in\mathbb{R}^n$, define $b_i(x):=\sum_{j=1}^n\textbf{1}_{\{x_i=x_{(j)}\}}\delta_j$ and $\varphi_i(x):=\sum_{j=1}^n\textbf{1}_{\{x_i=x_{(j)}\}}\sigma_j$, where $x_{(1)}\geq x_{(2)}\geq\cdots\geq x_{(n)}$ are ordered configuration of $\{x_1,x_2,\cdots,x_n\}$ with ties resolved by resorting to the lowest order.
Let $\phi$ be a nonnegative smooth function   on $\mathbb{R}^{n}$ supported in the unit ball with $\int_{\mathbb{R}^n}\phi(x)dx=1$.
Define
$$
b_i^m(x):=m^n\int_{\mathbb{R}^n}b_i(x-y)\phi(my)dy \quad \hbox{and} \quad
\varphi_i^{m}(x):=m^n\int_{\mathbb{R}^n}\varphi_i(x-y)\phi(my)dy.
$$
It is clear that $b_i^m$ and $\varphi_i^m$ are uniformly bounded smooth functions
with bounded first derivatives, and
$$\lim_{m\rightarrow \infty}b_i^m=b_i,\quad \lim_{m\rightarrow \infty}\varphi_i^m=\varphi_i
\quad \hbox{a.e. on } \mathbb{R}^n.
$$
 Moreover, $\varphi_i^m$ is bounded from below by a positive constant that is independent of
 $m\geq 1$.
It follows from \cite{BP} that uniqueness in law holds for SDE \eqref{1} and \cite{IKS} constructed a strong solution for SDE \eqref{1}. Consequently
(see, e.g., \cite[Theorem 3.2]{Cherney}), pathwise uniqueness holds for SDE \eqref{1}. Therefore, by Lemma \ref{lemma1}, we have for every $i=1,\ldots,n$, $p\geq1$,
and any compact subset $\mathcal{K}\subset \Gamma^n$,
$$\lim_{m\rightarrow \infty}\sup_{0\leq t\leq T}\sup_{x\in\mathcal{K}}\mathbb{E}\bigg[\sup_{0\leq s\leq T}\big|X_i^{t,x,m}(s)-X_i^{t,x}(s)\big|^p\bigg]=0,
$$
where
$X_i^{t,x,m}$ is the unique solution of
$$X_i^{t,x,m}(s)=x_i+\int_t^{s\vee t}b_i^m\big(X_i^{t,x,m}(u)\big)du+\int_t^{s\vee t}\varphi_i^m\big(X_i^{t,x,m}(u)\big)dW_i(u).$$
For any $x\in\Gamma^n$, since $\lim_{k\rightarrow\infty}|{x^{(k)}}-x|=0$, we can find a compact $\mathcal K\in\Gamma^n$ such that $x\in\mathcal K$ and ${x^{(k)}}\in\mathcal K$ for all $k\in\mathbb N$.
Therefore, for any $\varepsilon>0$, there exists $M\in\mathbb{N}$ (independent of $k$) such that for all $m\geq M$,
$$
\mathbb{E}\bigg[\sup_{0\leq s\leq T}\big|X_i^{t,x,m}(s)-X_i^{t,x}(s)\big|^p\bigg]\leq\frac{\varepsilon}{2},$$and $$\mathbb{E}\bigg[\sup_{0\leq s\leq T}\big|X_i^{{t^{(k)}},{x^{(k)}},m}(s)-X_i^{{t^{(k)}},{x^{(k)}}}(s)\big|^p\bigg]\leq\frac{\varepsilon}{2}.
$$
Note that $b_i^M$ and $\varphi_i^M$ are Lipschitz continuous, thus
\begin{equation}\nonumber
\lim_{k\rightarrow\infty}\mathbb{E}\bigg[\sup_{0\leq s\leq T}\big|X_i^{t,x,M}(s)-X_i^{{t^{(k)}},{x^{(k)}},M}(s)\big|^p\bigg]=0.
\end{equation}
Since $\varepsilon>0$ is arbitrary and
\begin{eqnarray*}
& &\mathbb{E}\bigg[\sup_{0\leq s\leq T}\big|X_i^{t,x}(s)-X_i^{{t^{(k)}},{x^{(k)}}}(s)\big|^p\bigg]\\
&\leq&C\mathbb{E}\bigg[\sup_{0\leq s\leq T}\big|X_i^{t,x}(s)-X_i^{t,x,M}(s)\big|^p\bigg]+C\mathbb{E}\bigg[\sup_{0\leq s\leq T}\big|X_i^{t,x,M}(s)-X_i^{{t^{(k)}},{x^{(k)}},M}(s)\big|^p\bigg]\\
& &+C\mathbb{E}\bigg[\sup_{0\leq s\leq T}\big|X_i^{{t^{(k)}},{x^{(k)}}}(s)-X_i^{{t^{(k)}},{x^{(k)}},M}(s)\big|^p\bigg]\\
&\leq&C\mathbb{E}\bigg[\sup_{0\leq s\leq T}\big|X_i^{t,x,M}(s)-X_i^{{t^{(k)}},{x^{(k)}},M}(s)\big|^p\bigg]+\varepsilon,
\end{eqnarray*}
we have \begin{equation}\nonumber
\lim_{k\rightarrow\infty}\mathbb{E}\bigg[\sup_{0\leq s\leq T}\big|X_i^{t,x}(s)-X_i^{{t^{(k)}},{x^{(k)}}}(s)\big|^p\bigg]=0.
\end{equation}
This proves the theorem.
\qed

By Theorem \ref{th1} and \cite[Lemma 2.5]{CF}, we can obtain the similar properties for ranked particles.
\begin{theorem}\label{th2}
Under the conditions of Theorem \ref{T:2.2}, for all $p\geq1$, there exists a constant depending on $\big(p,T,n,\{\delta_j\},\{\sigma_j\}\big)$ such that  for any $x,\{{x^{(k)}}\}_{k=1}^\infty\in\Gamma^n$ and for any $t,\{{t^{(k)}}\}_{k=1}^\infty\in[0,T]$ with $\lim_{k\rightarrow\infty}(|{t^{(k)}}-t|+|{x^{(k)}}-x|)=0$, we have for $j=1,\cdots,n,$
\begin{equation}\label{}
\mathbb{E}\bigg[\sup_{0\leq s\leq T}\big|X_{(j)}^{t,x}(s)\big|^p\bigg]\leq C\big(1+|x|^p\big),
\end{equation}
and
\begin{equation}\label{8}
\lim_{k\rightarrow\infty}\mathbb{E}\bigg[\sup_{0\leq s\leq T}\big|X_{(j)}^{t,x}(s)-X_{(j)}^{{t^{(k)}},{x^{(k)}}}(s)\big|^p\bigg]=0.
 \end{equation}
\end{theorem}

\subsection{BSDE with Rank-based Data}\label{S:2.3}

Throughout the rest of this paper, we always assume that,
the sequence $\{\sigma_1^2,
\dots, \sigma_n^2\}$ is concave. So by Theorem \ref{th1}, SDE \eqref{1}
with rank-based drift and diffusion coefficients has a unique strong solution
$ {X}^{t,x}(s):=\big(X_1^{t,x}(s),\cdots,X_n^{t,x}(s)\big)$ and that, with probability one,
there are no triple collisions at any time.
 Denote the ranked particles by
$$
\widetilde{X}^{t,x}(s):=\big(X_{(1)}^{t,x}(s),\cdots,X_{(n)}^{t,x}(s)\big).
$$

For any initial data $(t,x)\in[0,T]\times\Gamma^n$, consider the following BSDE (with $s$ running from $t$ to $T$):
\begin{equation}\label{e:2.11}
Y^{t,x}(s)=g\big(\widetilde{X}^{t,x}(T)\big)+\int_s^T
F \big(r,\wt X^{t,x}(r),Y^{t,x}(r), \bar Z^{t,x}(r)\big)dr-\int_s^T Z^{t,x}(r)\cdot dW(r) ,
\end{equation}
where $\bar Z^{t, {x}}_j(r):=\sum_{i=1}^n Z^{t, {x}}_i(r)1_{\{ X^{t, {x}}_i(r)=\widetilde X^{t, {x}}_j(r)\}}$,
and $F:[0,T]\times \Gamma^n \times\mathbb{R}\times\mathbb{R}^n\rightarrow\mathbb{R}$
and  $g:\Gamma^n\rightarrow\mathbb{R}$ are two  functions satisfying conditions ({\bf H1})-({\bf H2}) stated in Section 1.

Let $\sX$ be the set of all points in ${\mathbb R}^n$ that do not have more than two identical  coordinates.
Define for $(t, x, y, z)\in [0, \infty) \times \sX \times {\mathbb R} \times {\mathbb R}^n$,
$$
f(t, x, y, z)= F(t, \wt x, y, \bar z),
$$
where $\widetilde{x}=(x_{(1)},\ldots,x_{(n)})\in\Gamma^n$ is the rearrangement of coordinates of
$x$ in decreasing order,
and $\bar{z}_j:=\sum_{i=1}^nz_i\textbf{1}_{\{x_i=x_{(j)}\}} $.
Under condition (\textbf{H1}), $f$ is a jointly continuous function on $[0, \infty) \times \sX \times {\mathbb R} \times {\mathbb R}^n$, and has properties \eqref{12}-\eqref{13}.
We can rewrite BSDE \eqref{e:2.11} as
\begin{equation}\label{10}
Y^{t,x}(s)=g\big(\widetilde{X}^{t,x}(T)\big)+\int_s^Tf\big(r,X^{t,x}(r),Y^{t,x}(r),Z^{t,x}(r)\big)dr-\int_s^T Z^{t,x}(r)\cdot dW(r) ,
\end{equation}
Since the diffusion process $X^{t, x}$ takes values in $\sX$, we have by the same proof as that of   \cite[Theorem 4.1]{PP1} that  the BSDE \eqref{10}, and hence \eqref{e:2.11}, has a unique solution.
Recall the definition of the process $\{\beta_j (t), j=1, \dots, n\}$ from \eqref{e:2.3}, which are independent one-dimensional Brownian motions. Observe that
\begin{eqnarray*}
\int_s^T Z^{t,x}(r)\cdot dW(r)&=&\sum_{i=1}^n\int_s^TZ_i^{t,x}(r)dW_i(r)\\
&=&\sum_{i=1}^n\int_s^TZ_i^{t,x}(r)\big(\sum_{j=1}^n\textbf{1}_{\{X_i^{t,x}(r)=X_{(j)}^{t,x}(r)\}}d\beta_j(r)\big)\\
&=&\sum_{j=1}^n\int_s^T\sum_{i=1}^nZ_i^{t,x}(r)\textbf{1}_{\{X_i^{t,x}(r)=X_{(j)}^{t,x}(r)\}}d\beta_j(r)\\
&=&\sum_{j=1}^n\int_s^T\bar{Z}_j^{t,x}(r)d\beta_j(r).
\end{eqnarray*}
Thus BSDE \eqref{e:2.11} can be rewritten as
\begin{equation}\label{e:2.16}
Y^{t,x}(s)=g\big(\widetilde{X}^{t,x}(T)\big)+\int_s^TF\big(r,\widetilde{X}^{t,x}(r),Y^{t,x}(r),\bar{Z}^{t,x}(r)\big)dr-\int_s^T \bar{Z}^{t,x}(r)\cdot d\beta(r).
\end{equation}

It follows from \cite[Theorem 3.3]{CF} that BSDE \eqref{e:2.16} has a unique solution $(Y,\bar{Z})\in S^2\big([t,T];\mathbb{R}\big)\times M^2\big([t,T];\mathbb{R}^n\big)$.
We now define
\begin{equation}\label{16}
u(t,x):=Y^{t,x}(t),\ \ (t,x)\in[0,T]\times\Gamma^n,
\end{equation}
which is a deterministic quantity. In view of Theorem \ref{th1} and Theorem \ref{th2}, we can use the same argument as that of \cite[Theorem 3.5]{CF} to show that $(s,t,x)\rightarrow Y^{t,x}(s)$ is mean-square continuous. Consequently, $u(t,x)\in C\big([0,T]\times\Gamma^n\big).$
By an argument similar to that for \cite[Theorem 4.2 and 4.3]{CF}, we get the following result.

\begin{theorem}\label{th5}
Suppose {\bf(H1)} and {\bf(H2)} hold, $u(t,x)$ defined by \eqref{16} is a viscosity solution of the following PDE:
\begin{equation}\label{17}
\left\{
\begin{aligned}
& \frac{\partial u}{\partial t}(t,x) + \mathcal{L} u(t,x)
+ F\big(t,x,u(t,x),\sigma(\nabla u)(t,x)\big)=0
\quad \hbox{for } t\in[0,T) \hbox{ and } x\in \Pi^n, \\
&u(T,x)=g(x) \quad \hbox{for } x\in \Gamma^n,\\
& \frac{\partial u}{\partial x_{i+1}}(t,x)=\frac{\partial u}{\partial x_i}(t,x) \quad t\in[0,T) \hbox{ and } x\in F_i,\ i=1,\ldots,n-1.\\
\end{aligned}
\right.
\end{equation}
where
\begin{equation}\label{18}
\mathcal{L}=\frac{1}{2}\sum_{i=1}^n\sigma_i^2\frac{\partial^2}{\partial x_{i}^2}+\sum_{i=1}^n\delta_i\frac{\partial}{\partial x_{i}}
\end{equation}
and $\sigma$ is a diagonal matrix with diagonal elements $\sigma_1,\ldots,\sigma_n$. Furthermore, if for each $R>0$, there exists a positive function $\eta_R(\cdot)$ on $[0,\infty)$ with $\lim_{r\rightarrow0^+}\eta_R(r)=0$ such that
$$\big|F(t,x,y,z)-F(t,x',y,z)\big|\leq\eta_R\big(|x-x'|(1+|z|)\big),$$
for all $t\in[0,T]$, $z\in \mathbb{R}^n$, $x, x'\in \Gamma^n$ and $y\in \mathbb{R}$ having
$\max\{|x|,|x'|,|y|\}\leq R$, then $u(t,x)$ is the unique viscosity solution of \eqref{17} such that
\begin{equation}\label{19}
\lim_{|x|\rightarrow \infty}\big|u(t,x)\big|e^{-A\log^2|x|}=0,
\end{equation}
 uniformly for $t\in[0,T]$, for some $A>0$.
\end{theorem}

\section{ Reflected BSDE with Rank-based Data}\label{S:3}

For each $(t,x)\in[0,T]\times\Gamma^n$,
let $X^{t,x}(\cdot)$ be the solution to the SDE \eqref{1} with rank-based drifts $\{\delta_j\}$ and diffusion coefficients $\{\sigma_j\}$
and $\wt X^{t,x}(\cdot)$ the ranked process obtained from $X^{t,x}(\cdot)$ (see \S \ref{S:2.3}).
In this section, we study  reflected BSDE \eqref{e:1.1}.

\medskip

With the preparation in Section \ref{S:2},  we can use the arguments in
\cite[Theorem 5.2]{EKP} to establish the following result.

\begin{theorem}\label{th6}
Under {\bf(H1)}-{\bf(H3)}, there exists a unique triple $(Y^{t,x},Z^{t,x},K^{t,x})$ of progressively measurable processes such that
\begin{description}
  \item[(i)] $\mathbb{E}\sup_{t\leq s\leq T}|Y^{t,x}(s)|^2+\mathbb{E}\int_t^T|Z^{t,x}(s)|^2ds<\infty;$
  \item[(ii)] $Y^{t,x}(s)\geq h\big(s,\widetilde{X}^{t,x}(s)\big),$ $t\leq s\leq T;$
  \item[(iii)] $K^{t,x}(s)$ is a continuous increasing process and $$\int_t^T\big[Y^{t,x}(s)-h\big(s,\widetilde{X}^{t,x}(s)\big)\big]dK^{t,x}(s)=0.$$
  \item[(iv)] $(Y^{t,x},Z^{t,x},K^{t,x})$ satisfy \eqref{e:1.1}.
\end{description}
\end{theorem}

\medskip
Similar to \eqref{e:2.16}, \eqref{e:1.1} can be rewritten as
\begin{equation}\label{e:3.4}
\begin{aligned}
Y^{t,x}(s)&=g\big(\widetilde{X}^{t,x}(T)\big)+\int_s^TF\big(r,\widetilde{X}^{t,x}(r),Y^{t,x}(r),\bar{Z}^{t,x}(r)\big)dr+K^{t,x}(T)-K^{t,x}(s)\\
&\ \ \ -\int_s^T  \bar{Z}^{t,x}(r)\cdot d\beta(r) .
\end{aligned}
\end{equation}
Since
$$\mathbb{E}\sup_{t\leq s\leq T}|Y^{t,x}(s)|^2+\mathbb{E}\int_t^T|\bar{Z}^{t,x}(s)|^2ds< \infty,$$
and (ii), (iii) hold, we obtain that the triple $(Y^{t,x},\bar{Z}^{t,x},K^{t,x})$ is a solution of reflected BSDE \eqref{e:3.4} with obstacle $L(s):=h\big(s,\widetilde{X}^{t,x}(s)\big)$. The uniqueness is derived from that of the solution of reflected BSDE \eqref{e:1.1}.
We now define
\begin{equation}\label{24}
u(t,x):=Y^{t,x}(t),\ \ \ (t,x)\in[0,T]\times\Gamma^n.
\end{equation}

\begin{theorem}\label{th7}
Suppose {\bf(H1)}-{\bf(H3)} hold, we have $u\in C\big([0,T]\times\Gamma^n\big).$
\end{theorem}
\pf
It follows from \cite[Proposition 3.6]{EKP}, Theorem \ref{th1}, Theorem \ref{th2} and the continuity properties of $g,f,h$ that
$$\mathbb{E}\bigg[\sup_{t\leq s\leq T}\big|Y^{t',x'}(s)-Y^{t,x}(s)\big|^2\bigg]\rightarrow 0,\ \ as\ (t',x')\rightarrow(t,x).$$
Since $Y^{t,x}(t)$ is deterministic, we have $(t,x)\rightarrow Y^{t,x}(t)$ is continuous. \qed

\section{Viscosity solution of obstacle problem}\label{S:4}

In this section, we show that the solution of the reflected BSDE with rank-based data studied in the previous section can be used to give a probabilistic representation of the viscosity solution of the obstacle problem (or variational inequality) \eqref{e:1.7},
where  $\sigma_1,\cdots,\sigma_n$ and $\delta_1,\cdots,\delta_n$ are the diffusion and drift coefficients from \eqref{1}.

\medskip

First, we recall the definition of viscosity solution from  Crandall, Ishii and Lions \cite{CIL}.
Let $\delta=(\delta_1,\ldots,\delta_n)$ and $\sigma$ be the diagonal matrix with diagonal elements $\sigma_1,\ldots,\sigma_n$.

\begin{definition} \rm \begin{description}
\item{(i)} A function $u\in C\big([0,T]\times\Gamma^n\big)$ is called a \emph{viscosity subsolution} of \eqref{e:1.7} if
$$
u(T,x)\leq g(x)\quad \text{for}\ x\in\Gamma^n,
$$
and at any point $(t,x)\in[0,T)\times\Gamma^n$, for any $(p,Q,X)\in\bar{D}^{2,+}_u(t,x)$,
$$
\bigg(u(t,x)-h(t,x)\bigg)\wedge\bigg(-p-\frac{1}{2}
Tr(\sigma^2 X)
- \delta\cdot Q-F\big(t,x,u(t,x),\sigma Q\big)\bigg)\leq 0
\quad \text{for }\ x\in \Pi^n
$$
 and
$$
\bigg(Q_{i+1}-Q_i\bigg)\wedge\bigg\{\big(u(t,x)-h(t,x)\big)\wedge \big(-p-\frac{1}{2}Tr
(\sigma^2 X)
- \delta\cdot Q-F\big(t,x,u(t,x),\sigma Q\big)\big)\bigg\}\leq 0
$$
for $x\in F_i$,  $i=1,\ldots,n-1$.

\item{(ii)} A function $u\in C\big([0,T]\times\Gamma^n\big)$ is called a \emph{viscosity supersolution} of \eqref{e:1.7} if
$$
u(T,x)\geq g(x)\quad \text{for } x\in\Gamma^n,
$$
and at any point $(t,x)\in[0,T)\times\Gamma^n$, for any $(p,Q,X)\in\bar{D}^{2,-}_u(t,x)$,
$$\bigg(u(t,x)-h(t,x)\bigg)\wedge\bigg(-p-\frac{1}{2}
Tr(\sigma^2 X)
- \delta\cdot Q-F\big(t,x,u(t,x),\sigma Q\big)\bigg)\geq 0
\quad \text{for }  x\in \Pi^n
$$ and
$$\bigg(Q_{i+1}-Q_i\bigg)\vee\bigg\{\big(u(t,x)-h(t,x)\big)\wedge \big(-p-\frac{1}{2}Tr
( \sigma^2 X)
- \delta\cdot Q-F\big(t,x,u(t,x),\sigma Q\big)\big)\bigg\}\geq 0
$$
for $x\in F_i,$ $i=1,\ldots,n-1$.

\item{(iii)} A function $u\in C\big([0,T]\times\Gamma^n\big)$ is called a \emph{viscosity solution} of \eqref{e:1.7} if it is both a viscosity subsolution and a viscosity  supersolution.
\end{description}
\end{definition}

Next, we will give an equivalent definition of the viscosity solution.

\begin{definition} \rm \begin{description}
		\item{(i)} A function $u\in C\big([0,T]\times\Gamma^n\big)$ is called a \emph{viscosity subsolution} of  \eqref{e:1.7} if
		$$
		u(T,x)\leq g(x) \quad \hbox{for }  x\in\Gamma^n,
		$$
		and whenever $\varphi\in C^{1,2}\big([0,T]\times\Gamma^n\big)$ and  $(t,x)\in[0,T)\times\Gamma^n$ is a local minimum of $\varphi-u$, we have
		$$
		\bigg(u(t,x)-h(t,x)\bigg)\wedge\bigg(	 -\frac{\partial \varphi}{\partial t}(t,x)-\mathcal{L} \varphi(t,x)-F\big(t,x,u(t,x),\sigma\nabla \varphi(t,x)\big)\bigg)\leq0
		\quad \hbox{if } x\in \Pi^n
		$$
		and
		$$\bigg(\frac{\partial \varphi}{\partial x_{i+1}}(t,x)-\frac{\partial \varphi}{\partial x_i}(t,x)\bigg)\wedge\bigg\{\bigg(u(t,x)-h(t,x)\bigg)\wedge\bigg(- \frac{\partial \varphi}{\partial t}(t,x)-\mathcal{L} \varphi(t,x)-F\big(t,x,u(t,x),\sigma\nabla\varphi(t,x)\big)\bigg)\bigg\}\leq0
		$$
		if $x\in F_i$ for some  $i=1,\ldots,n-1$.
		
		\item{(ii)} A function $u\in C\big([0,T]\times\Gamma^n\big)$ is called a \emph{viscosity supersolution} of \eqref{e:1.7} if
		$$
		u(T,x)\geq g(x) \quad \hbox{for }  x\in\Gamma^n,
		$$
		and whenever $\varphi\in C^{1,2}\big([0,T]\times\Gamma^n\big)$ and  $(t,x)\in[0,T)\times\Gamma^n$ is a local maximum of $\varphi-u$, we have
		$$
		\bigg(u(t,x)-h(t,x)\bigg)\wedge\bigg(	 -\frac{\partial \varphi}{\partial t}(t,x)-\mathcal{L} \varphi(t,x)-F\big(t,x,u(t,x),\sigma\nabla \varphi(t,x)\big)\bigg)\geq0
		\quad \hbox{if } x\in \Pi^n
		$$
		and
		$$\bigg(\frac{\partial \varphi}{\partial x_{i+1}}(t,x)-\frac{\partial \varphi}{\partial x_i}(t,x)\bigg)\vee\bigg\{\bigg(u(t,x)-h(t,x)\bigg)\wedge\bigg(- \frac{\partial \varphi}{\partial t}(t,x)-\mathcal{L} \varphi(t,x)-F\big(t,x,u(t,x),\sigma\nabla\varphi(t,x)\big)\bigg)\bigg\}\geq0
		$$
		if $x\in F_i$ for some  $i=1,\ldots,n-1$.
		
		\item{(iii)} A function $u\in C\big([0,T]\times\Gamma^n\big)$ is called a \emph{viscosity solution} of \eqref{e:1.7} if it is both a viscosity subsolution and supersolution.
	\end{description}
\end{definition}

\begin{remark} \rm
(i) For $u\in C\big([0,T]\times\Gamma^n\big)$, if $u(t,x)\leq h(t,x)$ on $[0,T)\times\Gamma^n$ and $u(T,x)\leq g(x)$ on $\Gamma^n$, then clearly $u$ is a viscosity subsolution of \eqref{e:1.7}.

(ii) If a function $u\in C\big([0,T]\times\Gamma^n\big)$ is a viscosity supersolution of \eqref{e:1.7}, then $u(t,x) \geq h(t,x)$ on $[0,T)\times\Pi^n$ and $u(T,x)\geq g(x)$ on $\Gamma^n$.
\end{remark}

\subsection{Existence of viscosity solution}

In this subsection, we show that $u(t, x)$ defined by \eqref{24}
is a viscosity solution of the obstacle problem \eqref{e:1.7}.

\bigskip

\noindent {\bf Proof of Theorem \ref{th8}.}
First, note that $u(T,x)=g(x)$ for all $x\in\Gamma^n$. For each $(t,x)\in[0,T)\times\Gamma^n$, $m\in\mathbb{N}$, let $\big\{\big( Y^{t,x,m}(s),Z^{t,x,m}(s)\big),t\leq s\leq T\big\}$ denote the solution of BSDE
\begin{eqnarray*}
 Y^{t,x,m}(s)&=&g\big(\widetilde{X}^{t,x}(T)\big)+\int_s^TF\big(r,\widetilde{X}^{t,x}(r),Y^{t,x,m}(r),Z^{t,x,m}(r)\big)dr\\
& &+m\int_s^T\big( Y^{t,x,m}(r)-h(r,\widetilde{X}^{t,x}(r))\big)^{-}dr-\int_s^T Z^{t,x,m}(r)\cdot d\beta(r) ,\  t\leq s\leq T.\\
\end{eqnarray*}
It follows from Theorem \ref{th5} that $$u_m(t,x):= Y^{t,x,m}(t),\ \ 0\leq t\leq T,\ x\in\Gamma^n$$ is the viscosity solution of the parabolic PDE
\begin{equation}\nonumber
\left\{
\begin{aligned}
& \frac{\partial u_m}{\partial t}(t,x)
+ \mathcal{L} u_m(t,x) + F_m\big(t,x,u_m(t,x),\sigma(\nabla u_m)(t,x)\big)=0,\quad t\in[0,T),\ x\in \Pi^n, \\
&u_m(T,x)=g(x),\quad x\in \Gamma^n,\\
&\frac{\partial u_m}{\partial x_{i+1}}(t,x)
= \frac{\partial u_m}{\partial x_i}(t,x) \quad t\in[0,T),\ x\in F_i,\ i=1,\ldots,n-1,
\end{aligned}
\right.
\end{equation}
where
$$F_m(t,x,y,\sigma Q)=F(t,x,y,\sigma Q)+m\big(y-h(t,x)\big)^-.$$

\textbf{(i) Viscosity subsolution}: By the definition of viscosity subsolution, in order to prove $u$ is a viscosity subsolution, we only need to consider the case $u(t,x)> h(t,x)$ on $[0,T)\times\Gamma^n$. Define $\theta:=u(t,x)- h(t,x)>0$. Let $(p,Q,X)\in \bar{D}^{2,+}_u(t,x).$ By \cite[Lemma 6.1]{CIL}, there exist sequences$$ m_j\rightarrow \infty,\ (t^{(j)},x^{(j)})\rightarrow(t,x),\ (p^{(j)},Q^{(j)},X^{(j)})\in\bar{D}^{2,+}_{u_{m_j}}(t^{(j)},x^{(j)}),$$ such that $$(p^{(j)},Q^{(j)},X^{(j)})\rightarrow(p,Q,X).$$

If $x\in \Pi^n,$ we can find such sequence of $x^{(j)}$ such that $x^{(j)}\in \Pi^n$ for all $j$. Thus, for any $j$,
\begin{equation}\label{viscosity-1}
\begin{aligned}
&-p^{(j)}-\frac{1}{2}Tr(\sigma^2 X^{(j)})-\delta\cdot Q^{(j)}-F\big(t^{(j)},x^{(j)},u_{m_j}(t^{(j)},x^{(j)}),\sigma Q^{(j)}\big)\\&-m_j\big(u_{m_j}(t^{(j)},x^{(j)})-h(t^{(j)},x^{(j)})\big)^-\leq0.
\end{aligned}
\end{equation}

If there exists $i$ such that $1\leq i\leq n-1$ and $x\in F_i$, we can find an infinite subsequence of $\{x^{(j)}\}$ (still denoted by $\{x^{(j)}\}$), such that
$\{x^{(j)}\}\in \Pi^n$ for all $j$, or
$\{x^{(j)}\}\in F_i$ for all $j$.
If $\{x^{(j)}\}\in \Pi^n$ for all $j$, we have \eqref{viscosity-1}.
If $\{x^{(j)}\}\in F_i$ for all $j$,
then  for each $j$,
\begin{equation}\label{viscosity-2}
\begin{aligned}
\bigg(Q^{(j)}_{i+1}-Q^{(j)}_i\bigg)\wedge& \bigg\{-p^{(j)}-\frac{1}{2}Tr(\sigma^2 X^{(j)})-\delta\cdot Q^{(j)}-F\big(t^{(j)},x^{(j)},u_{m_j}(t^{(j)},x^{(j)}),\sigma Q^{(j)}\big)\\ &-m_j\big(u_{m_j}(t^{(j)},x^{(j)})-h(t^{(j)},x^{(j)})\big)^-\bigg\}\leq0.
\end{aligned}
\end{equation}
  Hence in either cases when $x\in F_i$ for some $1\leq i\leq n-1$, there is always an infinite subsequence of $\{x^{(j)}\}$
  (still denoted as $\{x^{(j)}\}$)
  so that \eqref{viscosity-2} holds for every $j\geq 1$.
It follows from Section 6 in \cite{EKP} that for each $0\leq s\leq T$, $y\in\Gamma^n$,
$$
u_m(s,y)\uparrow u(s,y)\quad  \hbox{as }  m\rightarrow \infty.
$$
Since $u_m$ and $u$ are continuous, it follows from Dini's theorem that the above convergence is uniform on compact set. Since $x^{(j)}\rightarrow x$, there exists a compact $\mathcal K$ such that $x,\{x^{(j)}\}_{j=1}^\infty\in\mathcal K$. It then exists a constant $M$ (independent of $(s,y)$) such that for $m\geq M$ and any $0\leq s\leq T$, $y\in\mathcal K$,
$$|u_m(s,y)-u(s,y)|\leq\theta /6.$$
It follows from the continuity of $u(\cdot,\cdot)$ and $h(\cdot,\cdot)$ that there exist $J_1$ and $J_2$ such that for $j\geq J_1$,
$$
|u(t,x)-u(t^{(j)},x^{(j)})|\leq\theta /6,
$$
 and for $j\geq J_2$,
$$
|h(t,x)-h(t^{(j)},x^{(j)})|\leq\theta /6.
$$ Since $m_j\rightarrow\infty$, there exists $J_3$ such that for $j\geq J_3$, $m_j\geq M$. Therefore, for $j\geq J:=J_1\vee J_2\vee J_3$,
 \begin{equation}\nonumber
\begin{aligned}
&u_{m_j}(t^{(j)},x^{(j)})-h(t^{(j)},x^{(j)})\\
=& \,  u_{m_j}(t^{(j)},x^{(j)})-u(t^{(j)},x^{(j)})+u(t^{(j)},x^{(j)})-u(t,x)+u(t,x)-h(t,x)+h(t,x)-h(t^{(j)},x^{(j)})\\
\geq& \,  {\theta}/{2}.
\end{aligned}
\end{equation}
 Therefore, taking the limit as $j\rightarrow \infty$, we obtain
for $x\in \Pi^n$,
$$-p-\frac{1}{2}Tr(\sigma^2 X)-\delta\cdot Q-F\big(t,x,u(t,x),\sigma Q\big)\leq0;$$
and for $x\in F_i$, $i=1,\ldots,n-1$,
$$
\bigg(Q_{i+1}-Q_i\bigg)\wedge \bigg(-p-\frac{1}{2}Tr(\sigma^2 X)-\delta\cdot Q-F\big(t,x,u(t,x),\sigma Q\big)\bigg)\leq0.$$

\textbf{(2) Viscosity supersolution}: Let $(t,x)$ be an arbitrary point in $[0,T)\times\Gamma^n$ and $(p,Q,X)\in\bar{D}^{2,-}_u(t,x)$. From the definition of the solution of reflected BSDE, we know that $u(t,x)\geq h(t,x)$ on $[0,T)\times\Gamma^n$.
By \cite[Lemma 6.1]{CIL}, there exist sequences
$$ m_j\rightarrow \infty,\quad (t^{(j)},x^{(j)})\rightarrow(t,x)
\quad \hbox{and} \quad  (p^{(j)},Q^{(j)},X^{(j)})\in\bar{D}^{2,-}_{u_{m_j}}(t^{(j)},x^{(j)})
$$
such that
$$
(p^{(j)},Q^{(j)},X^{(j)})\rightarrow(p,Q,X).
$$

If $x\in \Pi^n,$ we can find such a sequence of $x^{(j)}$ so that $x^{(j)}\in \Pi^n$ for all $j$. Thus, for any $j$,
\begin{equation}\label{viscosity-4}\begin{aligned}
	&
	-p^{(j)}-\frac{1}{2}Tr(\sigma^2 X^{(j)})-\delta\cdot Q^{(j)} -F\big(t^{(j)},x^{(j)},u_{m_j}(t^{(j)},x^{(j)}),\sigma Q^{(j)}\big)
	\\&  \quad
	 -m_j\big(u_{m_j}(t^{(j)},x^{(j)})-h(t^{(j)},x^{(j)})\big)^-\geq0.
\end{aligned}\end{equation}
Consequently,
\begin{eqnarray*}
	-p^{(j)}-\frac{1}{2}Tr(\sigma^2 X^{(j)})-\delta\cdot Q^{(j)} -F\big(t^{(j)},x^{(j)},u_{m_j}(t^{(j)},x^{(j)}),\sigma Q^{(j)}\big)
	\geq0.
\end{eqnarray*}

If there exists $i$ such that $1\leq i\leq n-1$ and $x\in F_i$, we can find an infinite subsequence of $\{x^{(j)}\}$ (still denoted by $\{x^{(j)}\}$), such that
$\{x^{(j)}\}\in \Pi^n$ for all $j$, or
$\{x^{(j)}\}\in F_i$ for all $j$.
If $\{x^{(j)}\}\in \Pi^n$ for all $j$, we have \eqref{viscosity-4}.
If $\{x^{(j)}\}\in F_i$ for all $j$, we have
\begin{equation}\label{viscosity-3}
\begin{aligned}
	\bigg(Q^{(j)}_{i+1}-Q^{(j)}_i\bigg)\vee& \bigg\{-p^{(j)}-\frac{1}{2}Tr(\sigma^2 X^{(j)})-\delta\cdot Q^{(j)}-F\big(t^{(j)},x^{(j)},u_{m_j}(t^{(j)},x^{(j)}),\sigma Q^{(j)}\big)\\& -m_j\big(u_{m_j}(t^{(j)},x^{(j)})-h(t^{(j)},x^{(j)})\big)^-\bigg\}\geq0.
\end{aligned}
\end{equation}
Therefore, in either case when $x\in F_i$ for some $1\leq i\leq n-1$, we always have \eqref{viscosity-3} for every $j$.
Consequently,
\begin{eqnarray*}
	\bigg(Q^{(j)}_{i+1}-Q^{(j)}_i\bigg)\vee \bigg\{-p^{(j)}-\frac{1}{2}Tr(\sigma^2 X^{(j)})-\delta\cdot Q^{(j)}-F\big(t^{(j)},x^{(j)},u_{m_j}(t^{(j)},x^{(j)}),\sigma Q^{(j)}\big)
	\bigg\}\geq0.
\end{eqnarray*}

Taking the limit as $j\rightarrow \infty$, we conclude that
if $x\in \Pi^n$,
$$-p-\frac{1}{2}Tr(\sigma^2 X)-\delta\cdot Q-F\big(t,x,u(t,x),\sigma Q\big)\geq0;$$
and if $x\in F_i$, $i=1,\ldots,n-1$,
$$
\bigg(Q_{i+1}-Q_i\bigg)\vee \bigg(-p-\frac{1}{2}Tr(\sigma^2 X)-\delta\cdot Q-F\big(t,x,u(t,x),\sigma Q\big)\bigg)\geq0.$$
This completes the proof. \qed

\subsection{Uniqueness of viscosity solution}

In this subsection, we present a proof of Theorem \ref{th9} on the uniqueness of viscosity solution.

\medskip

By definition, a viscosity solution $u$ of \eqref{e:1.7} is continuous and satisfies $u(t,x)\geq h(t,x)$ on $[0,T)\times\Pi^n$. For every $\alpha>0$, define
$$\Pi^{n,\alpha}:=\{x\in\Pi^n:d(x,\partial(\Pi^n))\geq\alpha\},$$
where $d(x,\partial(\Pi^n)):=\min_{y\in\partial(\Pi^n)}|x-y|.$ In order to prove Theorem \ref{th9}, it suffices to show that
there exists $\alpha_0>0$ so that for
any viscosity subsolution $u$ satisfying $u(t,x)\geq h(t,x)$ on $[0,T)\times\Pi^n$ and viscosity supersolution $v$
both satisfying
\eqref{28}, we have
$$
u(t,x)\leq v(t,x)\ \ \text{on}\ \ (0,T)\times \Pi^{n,\alpha}
\quad \hbox{for every } \alpha \in (0, \alpha_0).
$$
Denote the interior of $\Pi^{n,\alpha}$ by $int(\Pi^{n,\alpha})$.
For $i=1,\cdots,n-1$ and $k=1,\cdots,n-1$, define
$$B_i:=\left\{x=(x_1,\cdots,x_n)\in\Pi^{n,\alpha}\bigg|x_i=x_{i+1}+\sqrt{2}\alpha\right\},$$
and
$$A_k=\bigcup_{1\leq i_1<i_2<\cdots< i_k\leq n-1}\bigcap_{l=1}^kB_{i_l}.$$
Therefore, $$\Pi^{n,\alpha}=\bigg(\bigcup_{k=1}^{n-1}A_k\bigg)\bigcup\bigg(int(\Pi^{n,\alpha})\bigg).$$
We need the following two lemmas. Lemma \ref{lemma3} is proved in \cite{CF}.
\begin{lemma}\label{lemma3}
For every $ A>0$, there exists $C_1>0$ such that the function
$$
\Psi(t,x)=\exp\left( \big(C_1(T-t)+A\big)\psi(x)\right)
$$
satisfies
$$
-\frac{\partial\Psi}{\partial t}(t,x)-\mathcal{L}\Psi(t,x)-c\Psi(t,x)-c\big|\sigma D\Psi(t,x)\big|>0  \quad on\ [t_1,T]\times \Pi^{n},
$$
where $c$ is the Lipschitz constant of $F$,
$\psi(x)=\left(1+ \tfrac12 \log\big(|x|^2+1\big) \right)^2$
and $t_1 :=(T-A/C_1)^{+}$.
\end{lemma}

\begin{lemma}\label{lemma2}
Let $u$ be a viscosity subsolution and $v$ be a viscosity supersolution of \eqref{e:1.7} and $u(t,x)\geq h(t,x)$ on $[0,T)\times \Pi^n$.
Then for every $\alpha >0$,
the function $w:=u-v$ is a viscosity subsolution of
\begin{equation}\label{29}
\left\{
\begin{aligned}
&w(t,x)\wedge\bigg(-\frac{\partial w}{\partial t}(t,x)-\mathcal{L} w(t,x)-c\big(|w|+|\sigma\nabla w|\big)(t,x)\bigg)=0,\quad  (t,x)\in[0,T)\times {\rm int} (\Pi^{n,\alpha}),\\
&w(T,x)=0,\quad  x\in \Pi^{n,\alpha},\\
&\sum_{l=1}^k\left(\frac{\partial w}{\partial x_{i_l}}(t,x)-\frac{\partial w}{\partial x_{i_l+1}}(t,x)\right)=0,  \quad t\in[0,T),\  x\in A_k,\ k=1,\cdots,n-1,
\end{aligned}
\right.
\end{equation}
 i.e., $$
w(T,x)\leq 0\quad \text{for}\ x\in\Pi^{n,\alpha},
$$
and whenever $\varphi\in C^{1,2}\big([0,T]\times\Pi^{n,\alpha}\big)$ and  $(t,x)\in[0,T)\times\Pi^{n,\alpha}$ is a local minimum of $\varphi-w$, we have
		$$
		w(t,x)\wedge\bigg(	-\frac{\partial \varphi}{\partial t}(t,x)-\mathcal{L} \varphi(t,x)-c\big(|w|+|\sigma\nabla \varphi|\big)(t,x)\bigg)\leq0
		\quad \hbox{if } x\in int(\Pi^{n,\alpha})
		$$
		and
		$$\left(\sum_{l=1}^k\bigg(\frac{\partial \varphi}{\partial x_{i_l}}(t,x)-\frac{\partial \varphi}{\partial x_{i_l+1}}(t,x)\bigg)\right)\wedge\bigg\{w(t,x)\wedge\bigg(- \frac{\partial \varphi}{\partial t}(t,x)-\mathcal{L} \varphi(t,x)-c\big(|w|+|\sigma\nabla \varphi|\big)(t,x)\bigg)\bigg\}\leq0
		$$
		if $x\in A_k$, $k=1,\ldots,n-1$.
\end{lemma}
Before we give the proof of Lemma \ref{lemma2}, first we give a remark on the boundary condition of \eqref{29}.

\begin{remark}\rm
There are $C_{n-1}^k$ classes of points in $A_k$, i.e., there are $C_{n-1}^k$ combinations of $i_1,\cdots,i_k$. For example, if $i_l=l$, for $l=1,\cdots,k$, then
$$\sum_{l=1}^k\left(\frac{\partial w}{\partial x_{i_l}}(t,x)-\frac{\partial w}{\partial x_{i_l+1}}(t,x)\right)=\sum_{l=1}^k\left(\frac{\partial w}{\partial x_{l}}(t,x)-\frac{\partial w}{\partial x_{l+1}}(t,x)\right)=\frac{\partial w}{\partial x_{1}}(t,x)-\frac{\partial w}{\partial x_{k+1}}(t,x).$$
\end{remark}

\medskip

\noindent{\bf \textbf{Proof of Lemma \ref{lemma2}:}}
First, note that $$w(T,x)=u(T,x)-v(T,x)\leq0.$$
Let $\varphi\in C^{1,2}\big([0,T]\times\Pi^{n,\alpha}\big)$ and $(t_0,x_0)\in [0,T)\times\Pi^{n,\alpha}$ be a
local  minimum point of $\varphi-w$ with $w(t_0,x_0)=\varphi(t_0,x_0)$. Modifying $\varphi$ if necessary, we may assume without loss of generality that $(t_0,x_0)\in [0,T)\times\Pi^{n,\alpha}$ is a strict global minimum point of $\varphi-w$.
It suffices to show that $$w(t_0,x_0)\wedge\bigg(-\frac{\partial \varphi}{\partial t}(t_0,x_0)-\mathcal{L} \varphi(t_0,x_0)-c\big(|w|+|\sigma\nabla \varphi|\big)(t_0,x_0)\bigg)\leq0.$$
If $(t_0,x_0)$ is a point such that $u(t_0,x_0)=h(t_0,x_0)$, then $w(t_0,x_0)=u(t_0,x_0)-v(t_0,x_0)\leq0$. Therefore, we only need to consider the case $u(t_0,x_0)>h(t_0,x_0)$.
Define $$\psi_\varepsilon(t,x,y):=u(t,x)-v(t,y)-\frac{(x-y)^2}{\varepsilon^2}-\varphi(t,x),$$ where $\varepsilon$ is a positive parameter which tends to zero. Choose $R>0$ large enough and define $$\Pi^{n,\alpha,R}:=B_R\cap \Pi^{n,\alpha}$$ so that $(t_0,x_0)\in[0,T)\times\Pi^{n,\alpha,R}$, where $B_R$ is the open ball in $\mathbb{R}^n$ centered at origin with radius $R$. Let $(t_\varepsilon,x_\varepsilon,y_\varepsilon)$ be global maximum point of $\psi_\varepsilon(t,x,y)$ in $[0,T)\times\Pi^{n,\alpha,R}$, then by a classical argument in the theory of viscosity solution,
 we have
\begin{description}
  \item[(i)] $(t_\varepsilon,x_\varepsilon,y_\varepsilon)\rightarrow(t_0,x_0,x_0)
	\hbox{ as } \varepsilon\rightarrow0;$
  \item[(ii)] ${|x_\varepsilon-y_\varepsilon|^2}/{\varepsilon^2}$ is bounded and tends to zero as $\varepsilon\rightarrow0$.
\end{description}

Now, fix $\varepsilon>0$, it follows from \cite[Theorem 8.3]{CIL} that for any $\theta>0$ there exist $(X_\theta,Y_\theta)\in S(n)\times S(n)$ and $c_\theta\in\mathbb{R}$ such that
$$\big(c_\theta+\frac{\partial\varphi}{\partial t}(t_\varepsilon,x_\varepsilon),p_\varepsilon+\nabla\varphi(t_\varepsilon,x_\varepsilon),X_\theta\big)\in \bar{D}^{2,+}_u(t_\varepsilon,x_\varepsilon),$$
$$(c_\theta,p_\varepsilon,Y_\theta)\in\bar{D}^{2,-}_v(t_\varepsilon,y_\varepsilon),$$ and
$$\left(
    \begin{array}{cc}
      X_\theta & 0 \\
      0 & -Y_\theta \\
    \end{array}
  \right)\leq A+\theta A^2,$$ where $p_\varepsilon=\frac{2(x_\varepsilon-y_\varepsilon)}{\varepsilon^2}$,
$A=\left(
      \begin{array}{cc}
        D^2\varphi(t_\varepsilon,x_\varepsilon)+\frac{2}{\varepsilon^2} & -\frac{2}{\varepsilon^2} \\
        -\frac{2}{\varepsilon^2} & \frac{2}{\varepsilon^2} \\
      \end{array}
    \right).$

It is easy to check that
$$A+\theta A^2=\left(
                 \begin{array}{cc}
                   D^2\varphi(t_\varepsilon,x_\varepsilon) & 0 \\
                   0 & 0 \\
                 \end{array}
               \right)+\frac{2}{\varepsilon^2}\left(
                                                \begin{array}{cc}
                                                  I & -I \\
                                                  -I & I \\
                                                \end{array}
                                              \right)+\theta M(\varepsilon),$$ where
$$M(\varepsilon)=\left(
                  \begin{array}{cc}
                    (D^2\varphi(t_\varepsilon,x_\varepsilon))^2+\frac{4}{\varepsilon^2}D^2\varphi(t_\varepsilon,x_\varepsilon) & -\frac{2}{\varepsilon^2}D^2\varphi(t_\varepsilon,x_\varepsilon) \\
                    -\frac{2}{\varepsilon^2}D^2\varphi(t_\varepsilon,x_\varepsilon) & 0 \\
                  \end{array}
                \right)+\frac{8}{\varepsilon^4}\left(
                                                 \begin{array}{cc}
                                                   I & -I \\
                                                   -I & I \\
                                                 \end{array}
                                               \right).$$
Since $u$ and $v$ are viscosity subsolution and supersolution of \eqref{e:1.7}, respectively, we have
\begin{eqnarray*}
\bigg(u(t_\varepsilon,x_\varepsilon)-h(t_\varepsilon,x_\varepsilon)\bigg)&\wedge&\bigg(-c_\theta-\frac{\partial\varphi}{\partial t}(t_\varepsilon,x_\varepsilon)-\delta\cdot\big(p_\varepsilon+\nabla\varphi(t_\varepsilon,x_\varepsilon)\big)\\
& &-\frac{1}{2
}Tr(\sigma^2 X_\theta)-F\big(t_\varepsilon,x_\varepsilon,u(t_\varepsilon,x_\varepsilon),\sigma(p_\varepsilon+\nabla\varphi(t_\varepsilon,x_\varepsilon)\big)\big)\bigg)\leq0,
\end{eqnarray*}
and
$$\bigg(v(t_\varepsilon,x_\varepsilon)-h(t_\varepsilon,x_\varepsilon)\bigg)\wedge\bigg(-c_\theta-\delta\cdot p_\varepsilon-\frac{1}{2
}Tr(\sigma^2 Y_\theta)-F\big(t_\varepsilon,y_\varepsilon,v(t_\varepsilon,y_\varepsilon),\sigma p_\varepsilon\big)\bigg)\geq0.$$

Since $u(t_0,x_0)>h(t_0,x_0)$, $(t_\varepsilon,x_\varepsilon)\rightarrow(t_0,x_0)$, and $u,h$ are continuous, we have $u(t_\varepsilon,x_\varepsilon)>h(t_\varepsilon,x_\varepsilon)$ for sufficiently small $\varepsilon$. Note that $v(t,x)\geq h(t,x)$ on $[0,T)\times \Pi^n$, we have
$$-c_\theta-\frac{\partial\varphi}{\partial t}(t_\varepsilon,x_\varepsilon)-\delta\cdot\big(p_\varepsilon+\nabla\varphi(t_\varepsilon,x_\varepsilon)\big)-\frac{1}{2
}Tr(\sigma^2 X_\theta)-F\big(t_\varepsilon,x_\varepsilon,u(t_\varepsilon,x_\varepsilon),\sigma\big(p_\varepsilon+\nabla\varphi(t_\varepsilon,x_\varepsilon)\big)\big)\leq0,$$
and
$$-c_\theta-\delta\cdot p_\varepsilon-\frac{1}{2
}Tr(\sigma^2 Y_\theta)-F\big(t_\varepsilon,y_\varepsilon,v(t_\varepsilon,y_\varepsilon),\sigma p_\varepsilon\big)\geq0.$$
Thus,
\begin{eqnarray*}
0&\leq&\frac{\partial\varphi}{\partial t}(t_\varepsilon,x_\varepsilon)+\delta\cdot \nabla\varphi(t_\varepsilon,x_\varepsilon)
+\frac{1}{2}Tr(\sigma^2 X_\theta)-\frac{1}{2}Tr(\sigma^2 Y_\theta)\\& &+
F\big(t_\varepsilon,x_\varepsilon,u(t_\varepsilon,x_\varepsilon),\sigma\big(p_\varepsilon+\nabla\varphi(t_\varepsilon,x_\varepsilon)\big)\big)-
F\big(t_\varepsilon,y_\varepsilon,v(t_\varepsilon,y_\varepsilon),\sigma p_\varepsilon\big).
\end{eqnarray*}
First,
$$\frac{1}{2}Tr(\sigma^2 X_\theta)-\frac{1}{2}Tr(\sigma^2 Y_\theta)\leq\frac{1}{2}Tr\big(\sigma^2 D^2\varphi(t_\varepsilon,x_\varepsilon)\big)+\frac{\theta}{2
}R^\varepsilon(\sigma),$$where
$$R^\varepsilon(\sigma)=\bigg(\left(
                                \begin{array}{c}
                                  \sigma \\
                                  \sigma \\
                                \end{array}
                              \right),M(\varepsilon)\left(
                                                      \begin{array}{c}
                                                        \sigma \\
                                                        \sigma \\
                                                      \end{array}
                                                    \right)\bigg).$$
Finally,
\begin{eqnarray*}
& &F\big(t_\varepsilon,x_\varepsilon,u(t_\varepsilon,x_\varepsilon),\sigma\big(p_\varepsilon+\nabla\varphi(t_\varepsilon,x_\varepsilon)\big)\big)-
F\big(t_\varepsilon,y_\varepsilon,v(t_\varepsilon,y_\varepsilon),\sigma p_\varepsilon\big)\\
&=&F\big(t_\varepsilon,x_\varepsilon,u(t_\varepsilon,x_\varepsilon),\sigma\big(p_\varepsilon+\nabla\varphi(t_\varepsilon,x_\varepsilon)\big)\big)-
F\big(t_\varepsilon,x_\varepsilon,v(t_\varepsilon,y_\varepsilon),\sigma p_\varepsilon)\big)\\
& &+F\big(t_\varepsilon,x_\varepsilon,v(t_\varepsilon,y_\varepsilon), \sigma p_\varepsilon)\big)-F\big(t_\varepsilon,y_\varepsilon,v(t_\varepsilon,y_\varepsilon),\sigma p_\varepsilon\big)\\
&\leq&\eta\big(|x_\varepsilon-y_\varepsilon|(1+|\sigma p_\varepsilon|)\big)+c\big|u(t_\varepsilon,x_\varepsilon)-v(t_\varepsilon,y_\varepsilon)\big|+c\big| \sigma \nabla\varphi(t_\varepsilon,x_\varepsilon)\big|.\\
\end{eqnarray*}
We have denoted by $\eta$ the modulus $\eta_R$ which appears in \eqref{27} for $R$ large enough.
First letting $\theta\rightarrow0$ and then letting $\varepsilon\rightarrow0$, we obtain
$$-\frac{\partial\varphi}{\partial t}(t_0,x_0)-\delta\cdot \nabla\varphi(t_0,x_0)-\frac{1}{2}Tr\big(\sigma^2 D^2\varphi(t_0,x_0)\big)-c\big|w(t_0,x_0)\big|-c\big|\sigma \nabla\varphi(t_0,x_0)\big|\leq0.$$
Therefore, $$w(t_0,x_0)\wedge\bigg(-\frac{\partial \varphi}{\partial t}(t_0,x_0)-\mathcal{L} \varphi(t_0,x_0)-c\big(|w|+|\sigma\nabla \varphi|\big)(t_0,x_0)\bigg)\leq0.$$
\qed

\medskip
\noindent{\bf \textbf{Proof of Theorem \ref{th9}.}}
Suppose $u$ is a viscosity subsolution satisfying $u(t,x)\geq h(t,x)$ on $[0,T)\times\Pi^n$, $v$ is a viscosity supersolution and $u$, $v$ satisfy \eqref{28}.
Define $w:=u-v.$
From \eqref{28} we obtain that
$$
\lim_{|x|\rightarrow \infty}\big|w(t,x)\big| \,
e^{-A \left(1+\tfrac12 \log\left( |x|^2+1 \right) \right)^2}=0
$$
uniformly for $t\in[0,T]$ and for some $A>0$.
This implies that for every $\beta >0$ and $t\in [0, T]$,
 $$
\big|w(t,x)\big|<
\beta e^{A \left(1+\tfrac12 \log\left( |x|^2+1 \right) \right)^2}
\leq \beta\Psi(t,x)
\quad \hbox{when } |x| \hbox{ is large enough},
$$
where $\Psi (t,x)$ is the function defined in Lemma \ref{lemma3}.
Let $t_1 :=(T-A/C_1)^{+}$. Then
$$
M(\beta, T)
:=\max_{[t_1,T]\times \Pi^{n,\alpha}}(w-\beta\Psi)(t,x)e^{c(t-T)}
$$ is achieved at some point
$(\bar{t},\bar{x}) \in [t_1,T]\times   \Pi^{n,\alpha} $.
We claim that $M (\beta, T)\leq 0$ for any $\beta>0$.

If $\bar{t}=T$, since $$w(T,x)=u(T,x)-v(T,x)\leq0,$$ we have $$
w(t,x)-\beta\Psi(t,x)\leq0\ \text{on}\ [t_1,T]\times
\Pi^{n,\alpha},
$$
that is, $M (\beta, T) \leq0$ when $\bar{t}=T$.
We now consider the case  $\bar{t}<T$. Were $M (\beta, T)>0$, then
$$
w(\bar{t},\bar{x})={M(\beta, T)} e^{c(T-\bar{t})}+\beta \Psi(\bar{t},\bar{x})>0.
$$
Define
$$
\phi(t,x):=\beta\Psi(t,x)+(w-\beta\Psi)
(\bar{t},\bar{x})e^{c(\bar{t}-t)}.
$$
Clearly, $\phi (\bar{t}, \bar{x})=w(\bar{t}, \bar{x}) >0$.
Note that
$$
w(t,x)-\beta\Psi(t,x)\leq(w-\beta\Psi)(\bar{t},\bar{x})e^{c(\bar{t}-t)}.
$$
So $w-\phi$ attains a global maximum 0 in $[t_1,T)\times \Pi^{n,\alpha}$ at $(\bar{t},\bar{x})$.  Since $w$ is a viscosity subsolution of \eqref{29},  if $\bar{x}\in int(\Pi^{n,\alpha})$,
$$w(\bar{t},\bar{x})\wedge\left\{-\frac{\partial \phi}{\partial t}(\bar{t},\bar{x})-\mathcal{L} \phi(\bar{t},\bar{x})-c\big(w(\bar{t},\bar{x})+\big|\sigma\nabla \phi(\bar{t},\bar{x})\big|\big)\right\}\leq0;$$
if $\bar{x}\in A_k$ with some $1\leq k\leq n-1$,
$$\left(\sum_{l=1}^k\bigg(\frac{\partial \phi}{\partial x_{i_l}}(\bar{t},\bar{x})-\frac{\partial \phi}{\partial x_{i_l+1}}(\bar{t},\bar{x})\bigg)\right)\wedge w(\bar{t},\bar{x})\wedge\bigg(-\frac{\partial \phi}{\partial t}(\bar{t},\bar{x})-\mathcal{L} \phi(\bar{t},\bar{x})-c\big(w(\bar{t},\bar{x})+\big|\sigma\nabla \phi(\bar{t},\bar{x})\big|\big)\bigg)\leq0.$$
It is easy to see that for $i=1,\cdots,n$,
$$\frac{\partial \phi}{\partial x_i}(\bar{t},\bar{x})=2\beta(C_1(T-\bar{t})+A)\Psi(\bar{t},\bar{x})\sqrt{\psi(\bar{x})}\frac{\bar{x}_i}{1+|\bar{x}|^2}.$$
Then, for $\bar{x}\in A_k$ with some $1\leq k\leq n-1$, we have
$$\sum_{l=1}^k\bigg(\frac{\partial \phi}{\partial x_{i_l}}(\bar{t},\bar{x})-\frac{\partial \phi}{\partial x_{i_l+1}}(\bar{t},\bar{x})\bigg)>0.$$
Therefore, we have
 $$-\frac{\partial \phi}{\partial t}(\bar{t},\bar{x})-\mathcal{L} \phi(\bar{t},\bar{x})-c\big(w(\bar{t},\bar{x})+\big|\sigma\nabla \phi(\bar{t},\bar{x})\big|\big)\leq0.$$
The left-hand side of the above inequality is equal to
$$
\beta\bigg(-\frac{\partial \Psi}{\partial t}(\bar{t},\bar{x})-\mathcal{L} \Psi(\bar{t},\bar{x})-c\big(\Psi(\bar{t},\bar{x})+\big|\sigma\nabla \Psi(\bar{t},\bar{x})\big|\big)\bigg).
$$
This is impossible in view of Lemma \ref{lemma3} and so $M(\beta, T)\leq 0$ when $\bar{t}<T$.
This establishes the claim that $M\leq0$ for every $\beta>0$ and $0<\alpha\leq\alpha_0$.

Since $\beta>0$ is arbitrary, we conclude that
 $w(t,x)\leq 0$ on $[t_1,T]\times\Pi^{n,\alpha}$.
By the same argument but with $t_1$ in place of $T$, we have
$w(t,x)\leq 0$ on $[t_2,t_1]\times\Pi^{n}$ where
$t_2:=(t_1-A/C_1)^+$. Proceeding this in at most $C_1T/A$ number of steps,
we get $w(t,x)\leq 0$ on $(0,T)\times\Pi^{n,\alpha}$.
\qed

\begin{remark} \rm
It follows from
$$
\mathbb{E}\bigg(\sup_{t\leq s\leq T}\big|X_i^{t,x}(s)\big|^p\bigg)\leq C\big(1+|x|^p\big),\ i=1,\cdots,n,\ p\geq1,
$$
a standard estimation for reflected BSDE
(see \cite[p.735]{EKP})
$$
|Y^{t,x}(t)|^2\leq C\mathbb{E}\bigg[\big|g\big(\widetilde{X}^{t,x}(T)\big)\big|^2+\int_t^T\big|f\big(s,X^{t,x}(s),0,0\big)\big|^2 ds+\sup_{t\leq s\leq T}\big|h\big(s,\widetilde{X}^{t,x}(s)\big)\big|^2\bigg],
$$
and the continuity assumptions on $f,g$ and $h$
that
$u(t,x):=Y^{t,x}(t)$ grows at most polynomially at infinity. Therefore $u(t,x)$ is the unique viscosity solution of \eqref{e:1.7} in the class of viscosity solutions which satisfy \eqref{28} for some $A>0$.
\end{remark}

\section{American option pricing}\label{S:5}

In this section, we  study American option pricing problem with rank-based stock prices. First, define  $$(\Pi^n)^+:=\{x\in\Pi^n:x_1> x_2>\ldots> x_n>0\}.$$
Similarly, we define $(\Gamma^n)^+$, $(F_i)^+$, $i=1,\ldots,n-1$.

Consider a financial market $\mathcal{M}$ that consists of one bond and $n$ stocks. Fix $p=(p_0,p_1,\ldots, p_n)\in\mathbb{R}^+\times(\Gamma^{n})^+$ and $T>0$, let the prices $P^{t,p}_0(s), \{P^{t,p}_i(s)\}_{i=1}^n$ of these financial instruments evolve according to the following equations:
\begin{equation}\label{30}
\left\{
\begin{aligned}
&P^{t,p}_0(s)=p_0+\int_t^{s\vee t}P^{t,p}_0(u)r(u)du, \\
&P^{t,p}_i(s)=p_i+\int_t^{s\vee t}P^{t,p}_i(u)\Big(\sum_{j=1}^n1_{\{P^{t,p}_i(u)=P^{t,p}_{(j)}(u)\}}\delta_jdu+\sum_{j=1}^n1_{\{P^{t,p}_i(u)=P^{t,p}_{(j)}(u)\}}\sigma_jdW_i(u)\Big).
\end{aligned}
\right.
\end{equation}
Here, $r(s)$ (interest rate) is assumed to be a uniformly bounded deterministic function; $\delta_j,j=1,\ldots, n$ are real numbers and $\sigma_j,j=1,\ldots, n$ are positive real numbers and the sequence $\{\sigma_1^2,\ldots,\sigma_n^2\}$ is concave.
Define $X_i^{t,x}(s):=\log P_i^{t,p}(s)$, $1\leq i \leq n$. By It\^{o}'s formula,  for $1\leq i \leq n$,
\begin{equation}\label{31}
\left\{
\begin{aligned}
&dX_i^{t,x}(s)=\sum_{j=1}^n\textbf{1}_{\{X_i^{t,x}(s)=X_{(j)}^{t,x}(s)\}}(\delta_j-\frac{1}{2}\sigma_j^2)ds+\sum_{j=1}^n\textbf{1}_{\{X_i^{t,x}(s)=X_{(j)}^{t,x}(s)\}}\sigma_jdW_i(s),&s\geq t,\\
&X_i^{t,x}(s)=x_i=\log p_i,& s\leq t.
\end{aligned}
\right.
\end{equation}

Since $\{\sigma_1^2,\ldots,\sigma_n^2\}$ is concave, with probability one, there are no triple collisions of \eqref{31}. Therefore, as explained in \cite[Remark 2]{KPS}, \eqref{31} has a unique strong solution and
$$Leb\{s>t| \exists 1\leq i<j\leq n,X_i^{t,x}(s)=X_j^{t,x}(s)\}=0,$$
almost surely, where $Leb$ denotes the Lebesgue measure on $[0,\infty)$.
Therefore, it follows from \cite[Theorem 2.5]{BG} that the ranked log-price processes satisfy the following equations: for $s\geq t$,
\begin{equation}\nonumber
dX_{(j)}^{t,x}(s)=(\delta_j-\frac{1}{2}\sigma_j^2)ds+\sigma_j d\beta_j(s)+\frac{1}{2}\big(d\Lambda^{j,j+1}(s)-d\Lambda^{j-1,j}(s)\big),\ j=1,\ldots,n,
\end{equation}
where $\Lambda^{j,j+1}(s)$, $j=1,\cdots,n-1$ are the local times accumulated at the origin by the nonnegative semimartingales $$G_j(\cdot):=X_{(j)}^{t,x}(\cdot)-X_{(j+1)}^{t,x}(\cdot),\ \ j=1,\cdots,n-1,$$ over the time interval $[0,s]$, $\Lambda^{0,1}(\cdot)=\Lambda^{n,n+1}(\cdot)\equiv0,$ and
\begin{equation}\nonumber
\beta_j(\cdot):=\sum_{i=1}^n\int_0^{\cdot}\textbf{1}_{\{X_i^{t,x}(s)=X_{(j)}^{t,x}(s)\}}dW_i(s),\ j=1,\cdots,n.
\end{equation}
It follows from the existence and uniqueness of strong solution of \eqref{31} that \eqref{30} has a unique strong solution and
the ranked price processes satisfy the following equations:
\begin{equation}\nonumber
dP_{(j)}^{t,p}(s)=P_{(j)}^{t,p}(s)\big[\delta_jds+ \sigma_jd\beta_j(s)+\frac{1}{2}\big(d\Lambda^{j,j+1}(s)-d\Lambda^{j-1,j}(s)\big)\big],\ i=1,\ldots,n,
\end{equation}
Moreover, there exists a constant $C$ depending on $\big(T,\{\delta_i\},\{\sigma_i\},q\big)$ such that for $q\geq2$ and $i=1,\ldots,n$, $j=1,\ldots,n$,$$\mathbb{E}\bigg[\sup_{0\leq s\leq T}\big|P^{t,p}_i(s)|^q\bigg]\leq C\big(1+|p|^q\big), $$ and $$\mathbb{E}\bigg[\sup_{0\leq s\leq T}\big|P^{t,p}_{(j)}(s)|^q\bigg]\leq C\big(1+|p|^q\big). $$

We now consider the valuation problem of an American contingent claim $\{\xi(t); t\in[0,T]\}$, where the holder has the right to exercise the option at any time between 0 and $T$. It is well known that in general this claim cannot be hedged by a self-financing portfolio, and that it is necessary to introduce superhedging strategy.

\begin{definition}
A self-financing super-strategy is a vector process $(Y(\cdot),\pi(\cdot),C(\cdot)) $ where $Y(\cdot)$ is the market value (or wealth process), $\pi(\cdot)$ the portfolio process, and $C(\cdot)$ the cumulative consumption process, such that
$$ dY(t)=\Sigma^n_{i=1}\pi_i(t)\frac{dP_i(t)}{P_i(t)}+\big(Y(t)-\Sigma^n_{i=1}\pi_i(t)\big)\frac{dP_0(t)}{P_0(t)}-dC(t),\ \ \ \mathbb{E}\int_0^T\big|\pi(t)\big|^2dt< \infty.$$
$C(\cdot)$ is an increasing right-continuous, adapted process with $C(0)=0$.
Given a payoff $\{\xi(t);t\in[0,T]\}$, a super-strategy is called a superhedging strategy if
$$
Y(t)\geq\xi(t) \ \hbox{ for every }  t\in[0,T], \quad \mathbb{P}-a.s.$$
The smallest endowment to finance a superhedging strategy is the price of the American option.
\end{definition}
We refer the interested reader to El Karoui and  Quenez \cite{EQ} for more information on these terminologies.

\medskip

In our context, the wealth process satisfies the following equation:
\begin{equation}\label{32}
\left\{
\begin{aligned}
dY^{t,p}(s)&=\big[Y^{t,p}(s)r(s)+\sum_{i=1}^n\pi_i(s)\big(\sum_{j=1}^n\textbf{1}_{\{P_i^{t,p}(s)=P_{(j)}^{t,p}(s)\}}(\delta_j-r(s))\big)\big]ds-dC(s)\\
&\ \ \ \ +\sum_{i=1}^n\pi_i(s)\big(\sum_{j=1}^n\textbf{1}_{\{P_i^{t,p}(s)=P_{(j)}^{t,p}(s)\}}\sigma_j\big)dW_i(s), \\
Y^{t,p}(t)&=y,
\end{aligned}
\right.
\end{equation}
where constant $y>0$ represents the initial endowment.
From the definition of $\beta_j$, $1\leq j \leq n$, we have
\begin{equation}\nonumber
\begin{aligned}
dY^{t,p}(s)&=\left( Y^{t,p}(s)r(s)+\sum_{i=1}^n\pi_i(s)
\sum_{j=1}^n\textbf{1}_{\{P_i^{t,p}(s)=P_{(j)}^{t,p}(s)\}}(\delta_j-r(s)) \right)ds -d C(s)\\
&\ \ \ \ +\sum_{i=1}^n\pi_i(s)\left(\sum_{j=1}^n\textbf{1}_{\{P_i^{t,p}(s)=P_{(j)}^{t,p}(s)\}}\sigma_jd\beta_j(s)\right)\\
&=\left( Y^{t,p}(s)r(s)+\sum_{j=1}^n\big(\delta_j-r(s)\big) \sum_{i=1}^n\pi_i(s)\textbf{1}_{\{P_i^{t,p}(s)=P_{(j)}^{t,p}(s)\}} \right) ds-dC(s)\\
&\ \ \ \  +\sum_{j=1}^n\left(\sum_{i=1}^n\pi_i(s)\textbf{1}_{\{P_i^{t,p}(s)=P_{(j)}^{t,p}(s)\}}\right)\sigma_jd\beta_j(s).
\end{aligned}
\end{equation}
Define
$$
\bar{\pi}_j(s):=\sigma_j\sum_{i=1}^n\pi_i(s)\textbf{1}_{\{P_i^{t,p}(s)=P_{(j)}^{t,p}(s)\}},\ j=1,\ldots,n.
$$
Then
\begin{equation}\label{33}
dY^{t,p}(s)=\left( Y^{t,p}(s)r(s)+\sum_{j=1}^n\frac{\delta_j-r(s)}{\sigma_j}\bar{\pi}_j(s)\right) ds-dC(s)+\sum_{j=1}^n\bar{\pi}_j(s)d\beta_j(s).
\end{equation}
Consider the above equation with the obstacle $L(s):=h\big(s,P_0^{t,p}(s),\widetilde{P}^{t,p}(s)\big)$ and the terminal value $g\big(P_0^{t,p}(T),\widetilde{P}^{t,p}(T)\big)$, where $\widetilde{P}^{t,p}(s)=\big(P_{(1)}^{t,p}(s),\ldots,P_{(n)}^{t,p}(s)\big)$. The function $g:\mathbb{R}^+\times(\Gamma^n)^+\rightarrow\mathbb{R}^+$ satisfies
for all $ p_0,p_0'\in\mathbb{R}^+,$ $p,p'\in(\Gamma^n)^+$, there exists a constant $c$ such that
\begin{equation}\label{34}
\big|g(p_0,p)-g(p_0',p')\big|\leq c\big(|p_0-p_0'|+|p-p'|\big).
\end{equation}
The function $h:[0,T]\times \mathbb{R}^+\times(\Gamma^n)^+\rightarrow\mathbb{R}^+$ is jointly continuous and satisfies for all $ p_0\in\mathbb{R}^+,$ $p\in(\Gamma^n)^+$, there exist some $c>0, k\in\mathbb{N}$ such that
\begin{equation}\label{36}
h(t,p_0,p)\leq c\big(1+|p_0|^k+|p|^k\big),\ \ t\in[0,T].
\end{equation}
We assume moreover that for all $ p_0\in\mathbb{R}^+,$ $p\in(\Gamma^n)^+$,
\begin{equation}\label{37}
h(T,p_0,p)\leq g(p_0,p).
\end{equation}
From Section 3, we know that the following reflected BSDE
\begin{equation}\label{38}
\left\{
\begin{aligned}
&dY^{t,p}(s)=\big[Y^{t,p}(s)r(s)+\sum_{j=1}^n\frac{\delta_j-r(s)}{\sigma_j}\bar{\pi}_j(s)\big]ds-dC(s)+\sum_{j=1}^n\bar{\pi}_j(s)d\beta_j(s),\\
&Y^{t,p}(T)=g\big(P_0^{t,p}(T),\widetilde{P}^{t,p}(T)\big),
\end{aligned}
\right.
\end{equation}
admits a unique solution $(Y^{t,p},\bar{\pi},C)$. Therefore, for American contingent claim
\begin{equation}\label{39}
\xi(s)=\left\{
\begin{aligned}
&h\big(s,P_0^{t,p}(s),\widetilde{P}^{t,p}(s)\big),&\ t\leq s< T,\\
&g\big(P_0^{t,p}(T),\widetilde{P}^{t,p}(T)\big),&\ t=T,
\end{aligned}
\right.
\end{equation}
there exists a minimal square integrable superhedging strategy and $Y^{t,p}(s)$ is the price process of the contingent $\xi(s)$.

\section*{Acknowledgement}
The authors thank the anonymous referee for helpful comments that led to the improvement of  the exposition of this paper.

\vskip 0.3truein

\noindent {\bf Zhen-Qing Chen}

\smallskip \noindent
Department of Mathematics, University of Washington, Seattle,
WA 98195, USA.

\noindent E-mail: \texttt{zqchen@uw.edu}

\medskip

\noindent {\bf Xinwei Feng}

\smallskip \noindent
 Zhongtai Securities Institute for Financial Studies, Shandong University, Jinan, Shandong 250100, China.

 \noindent Email: {\texttt xwfeng@sdu.edu.cn}


\begin{thebibliography}{99}

\bibitem{BFK}A. D. Banner, R. Fernholz and I. Karatzas, Atlas models of equity markets. \emph{Ann. Appl. Probab.}, \textbf{15} (2005), 2296-2330.

\bibitem{BG}A. D. Banner and R. Ghomrasni, Local times of ranked continuous semimartingales. \emph{Stochastic Process Appl.}, \textbf{118} (2008), 1244-1253.

\bibitem{BP}R. F. Bass and E. Pardoux, Uniqueness for diffusions with piecewise constant coefficients. \emph{Probab. Theory Relat. Fields}, \textbf{76} (1987), 557-572.

\bibitem{CE}Z. Chen and L. Epstein, Ambiguity, risk, and asset returns in continuous time. \emph{Econometrica}, \textbf{70} (2002), 1403-1443.

\bibitem{CF}Z.-Q. Chen and X. Feng, Backward stochastic differential equations with rank-based data. {\it Sci. China Math. \bf 61} (2018), 27-56.

\bibitem{Cherney}A. S. Cherny, On the uniqueness in law and the pathwise uniqueness for stochastic differential equations. \emph{Theory Probab. Appl.}, \textbf{46} (2002), 406-419.

\bibitem{CIL}M. G. Crandall, H. Ishii and P. L. Lions, User's guide to viscosity solutions of second order partial differential equations. \emph{Bull. Amer. Math. Soc.}, \textbf{27} (1992), 1-67.

\bibitem{CK}J. Cvitani\'{c} and I. Karatzas, Backward stochastic differential equations with reflection and Dynkin games. \emph{Ann. Probab.}, \textbf{24} (1996), 2024-2056.

\bibitem{DW}J. Dai and R. J. Williams, Existence and uniqueness of semimartingale reflecting Brownian motions in convex polyhedrons. \emph{Theory Probab. Appl.}, \textbf{40} (1996), 1-40.

\bibitem{EKP}N. El Karoui, C. Kapoudjian, E. Pardoux, S. Peng and M. C. Quenez, Reflected solutions of backward SDE's, and related obstacle problems for PDE's. \emph{Ann. Probab.}, \textbf{25} (1997), 702-737.

\bibitem{EPQ}N. El Karoui, S. Peng and M. C. Quenez, Backward stochastic differential equations in finance. \emph{Math. Finance}, \textbf{7} (1997), 1-71.

\bibitem{EQ}N. El Karoui and M. C. Quenez, Non-linear pricing theory and backward stochastic differential equations. In \emph{Financial mathematics, Springer Berlin Heidelberg,} \textbf{1656} (1997), 191-246.

\bibitem{FIK}E. R. Fernholz, T. Ichiba, I. Karatzas and V. Prokaj, Planar diffusions with rank-based characteristics and perturbed Tanaka equations. \emph{Probab. Theory Relat. Fields}, \textbf{156} (2013), 343-374.

\bibitem{GM}I. Gyongy and T. Martnez, On stochastic differential equations with locally unbounded drift. \emph{Czech. Math. J.}, \textbf{51} (2001), 763-783.

\bibitem{HR}J. M. Harrison and M. I. Reiman, Reflected Brownian motion on an orthant. \emph{Ann. Probab.}, \textbf{9} (1981), 302-308.

\bibitem{IK}T. Ichiba and I. Karatzas, On collisions of Brownian particles. \emph{Ann. Appl. Probab.}, \textbf{20} (2010), 951-977.

\bibitem{IKS}T. Ichiba and I. Karatzas, Shkolnikov M. Strong solutions of stochastic equations with rank-based coefficients. \emph{Probab. Theory Relat. Fields}, \textbf{156} (2013), 229-248.

\bibitem{IW}N. Ikeda and S. Watanabe, Stochastic differential equations and diffusion processes. Elsevier, (1981).

\bibitem{KN}H. Kaneko and S. Nakao, A note on approximation for stochastic differential equations. \emph{S\'{e}minaire de Probabilit\'{e}s XXII. Springer Berlin Heidelberg}, \textbf{22} (1988), 155-162.

\bibitem{KPS}I. Karatzas, S. Pal and M. Shkolnikov, Systems of Brownian particles with asymmetric collisions, \emph{Ann. Inst. H. Poincar\'{e} Probab. Statist.}, \textbf{52} (2016), 323-354.

\bibitem{Krylov}N. V. Krylov, Controlled diffusion processes. Springer, (2008).

\bibitem{MC}J. Ma and J. Cvitani\'{c}, Reflected forward-backward SDE and obstacle problems with boundary conditions. \emph{J. Appl. Math. Stochastic Anal.}, \textbf{14} (2001), 113-138.

\bibitem{Pardoux}E. Pardoux, Backward stochastic differential equations and viscosity solutions of systems of semilinear parabolic and elliptic PDEs of second order. In \emph{Stochastic Analysis and Related Topics VI. Birkh\"{a}user Boston}, {\bf42} (1998), 79-127.

\bibitem{PP1}E. Pardoux and S. Peng, Adapted solution of a backward stochastic differential equation. \emph{Syst. Control Lett.}, \textbf{14} (1990), 55-61.

\bibitem{PP2}E. Pardoux and S. Peng, Backward stochastic differential equations and quasilinear parabolic partial differential equations. In \emph{Stochastic partial differential equations and their applications. Springer Berlin Heidelberg}, {\bf176} (1992), 200-217.

\bibitem{PZ}E. Pardoux and S. Zhang, Generalized BSDEs and nonlinear Neumann boundary value problems. \emph{Probab. Theory Relat. Fields}, \textbf{110} (1998),  535-558.

\bibitem{PENG}S. Peng, Backward stochastic differential equations and applications to optimal control. \emph{Appl. Math. Opt.}, \textbf{27} (1993), 125-144.

\bibitem{Sarantsev1}A. Sarantsev, Infinite systems of competing Brownian particles: existence, uniqueness and convergence results. arXiv preprint arXiv:1403.4229, (2014).

\bibitem{Sarantsev2}A. Sarantsev, Triple and simultaneous collisions of competing Brownian particles. \emph{Electron. J. Probab.}, \textbf{20} (2015), 1-28.

\bibitem{Shkolnikov}M. Shkolnikov, Competing particle systems evolving by interacting Levy processes. \emph{Ann. Appl. Probab.}, \textbf{21} (2011), 1911-1932.


\bibitem{Williams}R. J. Williams, Semimartingale reflecting Brownian motions in the orthant. \emph{IMA Vol. Math. Appl.}, \textbf{71} (1995), 125-137.

\end{thebibliography}
\end{document}